\documentclass{amsart}%
\usepackage{amsmath}
\usepackage{graphicx}
\usepackage{amscd}%
\usepackage{amsfonts}%
\usepackage{amssymb}
\newtheorem{theorem}{Theorem}[section]
\theoremstyle{plain}

\newtheorem{corollary}[theorem]{Corollary}

\newtheorem{example}[theorem]{Example}

\newtheorem{lemma}[theorem]{Lemma}

\newtheorem{problem}[theorem]{Question}
\newtheorem{proposition}[theorem]{Proposition}
\newtheorem{remark}[theorem]{Remark}

\numberwithin{equation}{section}

\usepackage[T1]{fontenc}

\begin{document}
\title[Gradient of a polynomial]{The gradient of a polynomial at infinity}
\author{Jacek Ch\k{a}dzy\'{n}ski}
\address{Faculty of Mathematics, University of \L\'{o}d\'{z}, 90-238 \L\'{o}d\'{z}, Poland}
\email{jachadzy@math.uni.lodz.pl}
\author{Tadeusz Krasi\'{n}ski}
\address{Faculty of Mathematics, University of \L\'{o}d\'{z}, 90-238 \L\'{o}d\'{z}, Poland}
\email{krasinsk@krysia.uni.lodz.pl}
\thanks{This research was partially supported by KBN Grant No 2 P03A 007 18}
\date{March 15, 2002}

\begin{abstract}
We give a full description of the growth of the gradient of a polynomial in
two complex variables at infinity near any fiber of the polynomial.

\end{abstract}
\maketitle

\section{Introduction}

Let $f:\mathbb{C}^{n}\rightarrow\mathbb{\mathbb{C}}$ be a non-constant
polynomial and let $\nabla f:\mathbb{C}^{n}\rightarrow\mathbb{\mathbb{C}}^{n}$
be its gradient. There exists a finite set $B(f)\subset\mathbb{\mathbb{C}}$
such that $f$ is a locally trivial $C^{\infty}$-bundle over
$\mathbb{\mathbb{C}}\setminus B(f).$ The set $B(f)$ is the union of the set of
critical values $C(f)$ of $f$ and critical values $\Lambda(f)$ corresponding
to the singularities of $f$ at infinity$.$ The set $\Lambda(f) $ is defined to
be the set of all $\lambda\in\mathbb{\mathbb{C}}$ for which there are no
neighbourhood $U$ of $\lambda$ and a compact set $K\subset\mathbb{C}^{n}$ such
that $f:f^{-1}(U)\setminus K\rightarrow U$ is a trivial $C^{\infty}$-bundle.
It is known that the set $\Lambda(f)$ is finite (\cite{Ph}, Appendix A1,
\cite{V}, Corollary 5.1). It turns out that for $\lambda\in\mathbb{\mathbb{C}%
}$ the property of being in $\Lambda(f)$ depends on the behaviour of the
gradient $\nabla f$ near the fiber $f^{-1}(\lambda).$

Ha in \cite{HA2} defined the notion of the \L ojasiewicz exponent
$\widetilde{\mathcal{L}}_{\infty,\lambda_{0}}(f)$ of the gradient $\nabla f$
at infinity near a fibre $f^{-1}(\lambda_{0})$ in the following way
\begin{equation}
\widetilde{\mathcal{L}}_{\infty,\lambda_{0}}(f):=\lim_{\delta\rightarrow0^{+}%
}\mathcal{L}_{\infty}(\nabla f|f^{-1}(D_{\delta})),\label{1.1'}%
\end{equation}
where $D_{\delta}:=\{\lambda\in\mathbb{\mathbb{C}}:\left|  \lambda-\lambda
_{0}\right|  <\delta\}$ and $\mathcal{L}_{\infty}(\nabla f|f^{-1}(D_{\delta
}))$ is the \L ojasiewicz exponent at infinity of the mapping $\nabla f$ on
the set $f^{-1}(D_{\delta})$ (see the definition in Section 3) and gave,
without a proof, a characterization of $\Lambda(f)$ for $n=2$ in terms of
$\widetilde{\mathcal{L}}_{\infty,\lambda_{0}}(f).$ Namely, $\lambda_{0}%
\in\Lambda(f)$ if and only if $\widetilde{\mathcal{L}}_{\infty,\lambda_{0}%
}(f)<0$ (or equivalently $\widetilde{\mathcal{L}}_{\infty,\lambda_{0}}%
(f)<-1$)$.$ A generalization of \ this result was given by Parusi\'{n}ski in
\cite{PA}. Moreover, Ha also gave a formula for $\widetilde{\mathcal{L}%
}_{\infty,\lambda_{0}}(f)$ in terms of Puiseux expansions of roots of the
polynomial $f-\lambda_{0}$ at infinity for $\lambda_{0}\in$ $\Lambda(f)$ (this
formula is analogous to the formula for the local \L ojasiewicz exponent of
the gradient $\nabla f$, given in \cite{KL}).

The aim of this paper is to give in the case $n=2$ a complete description of
the behaviour of the gradient $\nabla f$ near any fibre $f^{-1}(\lambda)$ for
$\lambda\in\mathbb{\mathbb{C}}.$ To achieve this we define a more convenient
\L ojasiewicz exponent at infinity of $\nabla f$ near a fibre $f^{-1}%
(\lambda)$ (it will turns out in Section 5 to be equivalent to the above one)
as the infimum of the \L ojasiewicz exponents at infinity of $\nabla f$ on
meromorphic curves ''approximating'' $f^{-1}(\lambda)$ at infinity. Precisely,
for a non-constant polynomial $f:\mathbb{C}^{n}\rightarrow\mathbb{\mathbb{C}}$
we define $\mathcal{L}_{\infty,\lambda}(f)$ by
\begin{equation}
\mathcal{L}_{\infty,\lambda}(f):=\inf_{\Phi}\frac{\deg\nabla f\circ\Phi}%
{\deg\Phi},\label{1.1}%
\end{equation}
where $\Phi=(\varphi_{1},...,\varphi_{n})$ is a meromorphic mapping defined in
a neighbourhood of $\infty$ in $\overline{\mathbb{\mathbb{C}}}$ such that
$\deg\Phi>0$ and $\deg(f-\lambda)\circ\Phi<0.$ Here $\deg\Phi:=\max
(\deg\varphi_{1},...,\deg\varphi_{n})$. We shall also call such mappings
meromorphic curves.

The main results of the paper are effective formulas for $\mathcal{L}%
_{\infty,\lambda}(f)$ for each $\lambda\in\mathbb{\mathbb{C}}$ and properties
of the function $\lambda\mapsto\mathcal{L}_{\infty,\lambda}(f) $ for $n=2.$ To
describe them we outline the contents of the sections.

Section 2 has an auxiliary character and contains technical results on
relations between roots of a polynomial and its derivatives.

In Section 3 we investigate $\mathcal{L}_{\infty,\lambda}(f)$ for $\lambda
\in\Lambda(f).$ In particular we obtain the all results of Ha with full proofs.

The main theorems are given in Section 4. They are Theorems \ref{T4.1},
\ref{T4.7} and \ref{T4.8} which give effective formulas for $\mathcal{L}%
_{\infty,\lambda}(f)$ for each $\lambda\in\mathbb{\mathbb{C}}$ in terms of the
resultant $\operatorname*{Res}_{y}(f(x,y)-\lambda,f_{y}^{\prime}(x,y)-u),$
where $\lambda,u$ are new variables and $(x,y)$ is a generic system of
coordinates in $\mathbb{C}^{2}$. As a consequence we obtain (Theorem
\ref{T4.9}) a basic property of the function $\lambda\mapsto\mathcal{L}%
_{\infty,\lambda}(f).$ Namely,
\begin{align*}
\mathcal{L}_{\infty,\lambda}(f)  & =const.\geqslant0\ \ \text{for}%
\ \lambda\notin\Lambda(f),\\
\mathcal{L}_{\infty,\lambda}(f)  & \in[-\infty,-1)\text{ \ \ \ for}%
\ \lambda\in\Lambda(f)
\end{align*}
Technical problems in proofs of Theorems \ref{T4.7} and \ref{T4.8} are caused
by the fact that the exponent $\mathcal{L}_{\infty,\lambda}(f)$ does not have
to be attained on a meromorphic curve i.e. it may happen that there is no
meromorphic curve $\Phi$ such that $\deg\Phi>0,$ $\deg(f-\lambda)\circ\Phi<0$
and
\[
\mathcal{L}_{\infty,\lambda}(f)=\frac{\deg\nabla f\circ\Phi}{\deg\Phi}.
\]
At the end of this section we compare the functions $\lambda\mapsto
\mathcal{L}_{\infty,\lambda}(f)$ and $\lambda\mapsto\mathcal{L}_{\infty
}(\nabla f|f^{-1}(\lambda))$ (Theorem $\ref{T6.1}).$ It turns out they differ
only in a finite set of points, containing $\Lambda(f)$.

In Section 5 we show the equivalence (only in case $n=2$) of the two above
notions of the \L ojasiewicz exponent at infinity of $\nabla f$ near a fibre
$f^{-1}(\lambda)$ i.e. \ the equality
\begin{equation}
\widetilde{\mathcal{L}}_{\infty,\lambda}(f)=\mathcal{L}_{\infty,\lambda
}(f)\;\;\;\;\text{for \ }\lambda\in\mathbb{\mathbb{C}}.\label{1.2}%
\end{equation}

In Section 6 some partial results on $\mathcal{L}_{\infty,\lambda}(f)$ and
$\widetilde{\mathcal{L}}_{\infty,\lambda}(f)$ in $n$ dimensional case are
given. Here a complete characterization (in terms of the exponents
$\mathcal{L}_{\infty,\lambda}(f)$ and $\widetilde{\mathcal{L}}_{\infty
,\lambda}(f)) $ of the values of $f$ for which the Malgrange condition does
not hold is given. Some open problems on $\widetilde{\mathcal{L}}%
_{\infty,\lambda}(f)$ and $\mathcal{L}_{\infty,\lambda}(f)$ are also posed.

In the end we explain some technical assumptions occured in Sections 2-5.
Since one can easily show that the exponent $\mathcal{L}_{\infty,\lambda}(f)$
does not depend on linear change of coordinates in $\mathbb{C}^{n}$ we shall
assume in Sections 2-5 that the polynomial $f\in\mathbb{C}[x,y]$ is monic with
respect to $y$ and $\deg f=\deg_{y}f.$ Then we have a simple characterization
of the set $\Lambda(f)$ which will be used in the paper. Namely, in \cite{HA1}
and \cite{K} there was proved that
\begin{equation}
\Lambda(f)=\{\lambda\in\mathbb{\mathbb{C}}:c_{0}(\lambda)=0\},\label{1.4}%
\end{equation}
where the polynomial $c_{0}(\lambda)x^{N}+\cdots+c_{N}(\lambda)$ is the
resultant of the polynomials $f(x,y)-\lambda,$ $f_{y}^{\prime}(x,y)$ with
respect to the variable $y.$

\section{Auxiliary results}

Let $f$ be a non-constant polynomial in two complex variables of the form
\begin{equation}
f(x,y)=y^{n}+a_{1}(x)y^{n-1}+\cdots+a_{n}(x),\;\;\deg a_{i}\leqslant
i,\;\;i=1,...,n.\label{2.1}%
\end{equation}

It can be easily showed (see \cite{CK1}).

\begin{lemma}
\label{L2.1}If $n>1,$ then for every $\lambda_{0}\in\mathbb{C}$ there exist an
$D\in\mathbb{N}$ and functions $\beta_{1},...,\beta_{n},$ $\gamma
_{1},...,\gamma_{n-1}$, meromorphic at infinity, such that
\end{lemma}

(a) $\ \deg\beta_{i}\leqslant D,\;\;\deg\gamma_{j}\leqslant D,$

(b) $\ f(t^{D},y)-\lambda_{0}=\prod_{i=1}^{n}(y-\beta_{i}(t)),$

(c) $\ f_{y}^{\prime}(t^{D},y)=n\prod_{j=1}^{n-1}(y-\gamma_{j}(t)).$

\begin{flushright}$\square$\end{flushright}

We shall now give a lemma which directly follows from the property B.3 in
\cite{GP}. Local version of this lemma was proved in \cite{KL}.

\begin{lemma}
\label{LD2}(\cite{GP}, B.3) Under notation and assumptions of Lemma \ref{L2.1}
for every $i,j\in\{1,...,n\},$ $i\neq j,$ there exists $k\in\{1,...,n-1\}$
such that
\begin{equation}
\deg(\beta_{i}-\beta_{j})=\deg(\beta_{i}-\gamma_{k})\label{D.1}%
\end{equation}
and inversly for every $i\in\{1,...,n\}$ and $k\in\{1,...,n-1\}$ there exists
$j\in\{1,...,n\}$ such that (\ref{D.1}) holds.
\end{lemma}

\begin{flushright}$\square$\end{flushright}

Now we prove a proposition useful in the sequel. A local version of it is
given in \cite{P2} and \cite{R}. We put $\Psi_{l}(t):=(t^{D},\gamma_{l}(t)),$
$l\in\{1,...,n-1\}.$

\begin{proposition}
\label{P2.2}Under notations and assumptions of Lemma \ref{L2.1} we have
\begin{equation}
\min_{i=1}^{n}\left(  \sum_{j=1,j\neq i}^{n}\deg\left(  \beta_{i}-\beta
_{j}\right)  +\min_{j=1,j\neq i}^{n}\deg(\beta_{i}-\beta_{j})\right)
=\min_{l=1}^{n-1}(\deg\left(  f-\lambda_{0}\right)  \circ\Psi_{l}).\label{2.2}%
\end{equation}
\end{proposition}

\begin{proof}
(after \cite{R}). There exists $i_{0}\in\{1,...,n\}$ such that the left hand
side in (\ref{2.2}) is equal to
\[
\sum_{j=1,j\neq i_{0}}^{n}\deg\left(  \beta_{i_{0}}-\beta_{j}\right)
+\min_{j=1,j\neq i_{0}}^{n}\deg(\beta_{i_{0}}-\beta_{j})
\]
and $j_{0}\in\{1,...,n\}$ such that
\begin{equation}
\min_{j=1,j\neq i_{0}}^{n}\deg(\beta_{i_{0}}-\beta_{j})=\deg(\beta_{i_{0}%
}-\beta_{j_{0}}).\label{2.3}%
\end{equation}
By Lemma \ref{LD2} there exists $k_{0}\in\{1,...,n-1\}$
\begin{equation}
\deg(\beta_{i_{0}}-\beta_{j_{0}})=\deg(\beta_{i_{0}}-\gamma_{k_{0}%
}).\label{2.4}%
\end{equation}
We shall lead the further part of the proof in four steps.

A. We first show that for each $j\in\{1,...,n\}$ we have
\begin{equation}
\deg(\gamma_{k_{0}}-\beta_{j})\geqslant\deg(\beta_{i_{0}}-\beta_{j_{0}%
}).\label{2.5}%
\end{equation}
Take any $j\in\{1,...,n\}$ and consider two cases:

(a) $\deg(\beta_{i_{0}}-\beta_{j_{0}})\leqslant\min_{\substack{s=1 \\s\neq
j}}^{n}\deg(\beta_{s}-\beta_{j}),$

(b) $\deg(\beta_{i_{0}}-\beta_{j_{0}})>\min_{\substack{s=1 \\s\neq j}}^{n}%
\deg(\beta_{s}-\beta_{j}).$

\noindent In case (a) by Lemma \ref{LD2} there exists $p\in\{1,...,n\}$ such
that
\[
\deg(\gamma_{k_{0}}-\beta_{j})=\deg(\beta_{p}-\beta_{j})\geqslant
\min_{\substack{s=1 \\s\neq j}}^{n}\deg(\beta_{s}-\beta_{j})\geqslant
\deg(\beta_{i_{0}}-\beta_{j_{0}}).
\]
which gives (\ref{2.5}).

In case (b) by definition of $i_{0}$ and (\ref{2.3}) we have
\begin{align*}
&  \sum_{s=1,s\neq i_{0}}^{n}\deg\left(  \beta_{s}-\beta_{i_{0}}\right)
+\deg(\beta_{j_{0}}-\beta_{i_{0}})\\
&  \leqslant\sum_{s=1,s\neq j}^{n}\deg(\beta_{s}-\beta_{j})+\min_{s=1,s\neq
j}^{n}\deg(\beta_{s}-\beta_{j}).
\end{align*}
Hence and from (b) we get
\[
\sum_{s=1,s\neq i_{0}}^{n}\deg\left(  \beta_{s}-\beta_{i_{0}}\right)
<\sum_{s=1,s\neq j}^{n}\deg(\beta_{s}-\beta_{j}).
\]
Then there exists $s\neq i_{0}$, $s\neq j$ such that
\[
\deg\left(  \beta_{s}-\beta_{i_{0}}\right)  <\deg(\beta_{s}-\beta_{j}).
\]
Hence and from (\ref{2.4}) we get
\begin{align*}
\deg(\beta_{j}-\beta_{i_{0}}) &  =\deg(\beta_{j}-\beta_{s}+\beta_{s}%
-\beta_{i_{0}})>\deg(\beta_{s}-\beta_{i_{0}})\\
&  \geqslant\deg(\beta_{j_{0}}-\beta_{i_{0}})=\deg(\gamma_{k_{0}}-\beta
_{i_{0}})
\end{align*}
In consequence
\[
\deg(\gamma_{k_{0}}-\beta_{j})=\deg(\gamma_{k_{0}}-\beta_{i_{0}}+\beta_{i_{0}%
}-\beta_{j})>\deg(\gamma_{k_{0}}-\beta_{i_{0}})=\deg(\beta_{j_{0}}%
-\beta_{i_{0}}).
\]
This gives (\ref{2.5}) in case (b).

B. We shall now show that for each $j\in\{1,...,n\}$, $j\neq i_{0},$ we have
\begin{equation}
\deg(\gamma_{k_{0}}-\beta_{j})=\deg(\beta_{i_{0}}-\beta_{j}).\label{2.6}%
\end{equation}
Take $j\in\{1,...,n\}$, $j\neq i_{0},$ and consider two cases:

(a) $\deg(\gamma_{k_{0}}-\beta_{j})>\deg(\gamma_{k_{0}}-\beta_{i_{0}}), $

(b) $\deg(\gamma_{k_{0}}-\beta_{j})=\deg(\gamma_{k_{0}}-\beta_{i_{0}}). $

By (\ref{2.4}) and (\ref{2.5}) there are no more cases. In case (a) we have
\[
\deg(\beta_{i_{0}}-\beta_{j})=\deg(\beta_{i_{0}}-\gamma_{k_{0}}+\gamma_{k_{0}%
}-\beta_{j})=\deg(\gamma_{k_{0}}-\beta_{j}),
\]
which gives (\ref{2.6}).

In case (b) by (\ref{2.3}) and (\ref{2.4}) we have
\begin{align*}
\deg(\beta_{i_{0}}-\beta_{j}) &  =\deg(\beta_{i_{0}}-\gamma_{k_{0}}%
+\gamma_{k_{0}}-\beta_{j})\leqslant\deg(\beta_{i_{0}}-\gamma_{k_{0}})\\
&  =\deg(\beta_{i_{0}}-\beta_{j_{0}})\leqslant\deg(\beta_{i_{0}}-\beta_{j}),
\end{align*}
which gives (\ref{2.6}) in case (b).

C. We notice that by (\ref{2.6}) we have
\begin{align*}
\deg(f-\lambda_{0})\circ\Psi_{k_{0}} &  =\sum_{j=1}^{n}\deg\left(
\gamma_{k_{0}}-\beta_{j}\right) \\
&  =\sum_{j=1,j\neq i_{0}}^{n}\deg\left(  \gamma_{k_{0}}-\beta_{j}\right)
+\deg\left(  \gamma_{k_{0}}-\beta_{i_{0}}\right) \\
&  =\sum_{j=1,j\neq i_{0}}^{n}\deg\left(  \beta_{i_{0}}-\beta_{j}\right)
+\min_{_{\substack{j=1 \\j\neq i_{0}}}}^{n}\deg\left(  \beta_{i_{0}}-\beta
_{j}\right)  .
\end{align*}
Thus we have shown
\begin{equation}
\min_{i=1}^{n}\left(  \sum_{j=1,j\neq i}^{n}\deg\left(  \beta_{i}-\beta
_{j}\right)  +\min_{j=1,j\neq i}^{n}\deg(\beta_{i}-\beta_{j})\right)
\geqslant\min_{l=1}^{n-1}(\deg\left(  f-\lambda_{0}\right)  \circ\Psi
_{l}).\label{2.7}%
\end{equation}

D. We shall now show the inequality opposite to (\ref{2.7}). There exist
\linebreak $l_{0}\in\{1,...,n-1\}$ and $j_{0}\in\{1,...,n\}$ such that
\begin{align}
\min_{l=1}^{n-1}(\deg\left(  f-\lambda_{0}\right)  \circ\Psi_{l}) &
=\deg\left(  f-\lambda_{0}\right)  \circ\Psi_{l_{0}},\label{2.8}\\
\min_{j=1}^{n}\deg(\gamma_{l_{0}}-\beta_{j}) &  =\deg(\gamma_{l_{0}}%
-\beta_{j_{0}}).\label{2.9}%
\end{align}

Observe first that for any $j\in\{1,...,n\},$ $j\neq j_{0},$ we have by
(\ref{2.9})
\begin{equation}
\deg\left(  \beta_{j}-\beta_{j_{0}}\right)  =\deg\left(  \beta_{j}%
-\gamma_{l_{0}}+\gamma_{l_{0}}-\beta_{j_{0}}\right)  \leqslant\deg(\beta
_{j}-\gamma_{l_{0}}).\label{2.10}%
\end{equation}
By Lemma \ref{LD2} there exists $k_{0}\in\{1,...,n-1\}$ such that $\deg
(\gamma_{l_{0}}-\beta_{j_{0}})=\deg\left(  \beta_{k_{0}}-\beta_{j_{0}}\right)
.$ Hence using Lemma \ref{L2.1} and (\ref{2.10}) we get
\begin{align*}
\deg\left(  f-\lambda_{0}\right)  \circ\Psi_{l_{0}} &  =\sum_{j=1}^{n}%
\deg\left(  \gamma_{l_{0}}-\beta_{j}\right)  \geqslant\sum_{j=1,j\neq j_{0}%
}^{n}\deg\left(  \beta_{j_{0}}-\beta_{j}\right)  +\deg(\gamma_{l_{0}}%
-\beta_{j_{0}})\\
&  \geqslant\sum_{j=1,j\neq j_{0}}^{n}\deg\left(  \beta_{j_{0}}-\beta
_{j}\right)  +\min_{j=1,j\neq j_{0}}^{n}\deg(\beta_{j}-\beta_{j_{0}}),
\end{align*}
which gives the inequality opposite to (\ref{2.7}).

This ends the proof.
\end{proof}

\section{Critical values at infinity}

We start with definitions and notation needed in the sequel.

Let $F:\mathbb{C}^{n}\rightarrow\mathbb{C}^{m}$, $n\geqslant2$, be a
polynomial mapping and let $S\subset\mathbb{C}^{n}$ be an unbounded set. We
define
\[
N(F|S):=\{\nu\in\mathbb{R}:\mathbb{\exists}A,B>0\;\forall z\in S\;(\left|
z\right|  >B\Rightarrow\left|  F(z)\right|  \geqslant A\left|  z\right|
^{\nu})\},
\]
where $\left|  \cdot\right|  $ is the polycylindric norm. If $S=\mathbb{C}^{n}
$ we put $N(F):=N(F|\mathbb{C}^{n}).$

By the \L ojasiewicz exponent at infinity of $F|S$ we shall mean
$\mathcal{L}_{\infty}(F|S):=\sup N(F|S)$ when $N(F|S)\neq\emptyset,$ and
$-\infty$ when $N(F|S)=\emptyset$. Analogously $\mathcal{L}_{\infty}(F):=\sup
N(F)$ when $N(F)\neq\emptyset,$ and $-\infty$ when $N(F)=\emptyset$.

We give now a lemma needed in the sequel, which gives known formulas for the
\L ojasiewicz exponent at infinity of a polynomial on the zero set of the
another one. Let $g,h$ be polynomials in two complex variables $(x,y)$ and
\[
0<\deg h=\deg_{y}h.
\]
Let $\tau\in\mathbb{C}$ and $R(x,\tau):=\operatorname*{Res}_{y}(g(x,y)-\tau
,h(x,y))$ be the resultant of $g(x,y)-\tau$ and $h(x,y)$ with respect to $y. $
We put
\begin{gather*}
R(x,\tau)=R_{0}(\tau)x^{K}+\cdots+R_{K}(\tau),\;\;\;R_{0}\neq0,\\
T:=h^{-1}(0).
\end{gather*}

\begin{lemma}
\label{LP}(\cite{P1}, Proposition 2.4) Under above notation and assumptions
there is:
\end{lemma}

\begin{description}
\item[(i)] $\mathcal{L}_{\infty}(g|T)>0$ \textit{if and only if} $R_{0}=const.,$

\item[(ii)] $\mathcal{L}_{\infty}(g|T)=0$ \textit{if and only if} $R_{0}\neq
const.$ \textit{and} $R_{0}(0)\neq0,$

\item[(iii)] $-\infty<\mathcal{L}_{\infty}(g|T)<0$ \textit{if and only if
there exists }$r$\textit{\ such that} $R_{0}(0)=\cdots=R_{r}(0)=0$
\textit{and} $R_{r+1}(0)\neq0,$

\item[(iv)] $L_{\infty}(g|T)=-\infty$\textit{\ if and only if }$R_{0}%
(0)=\cdots=R_{K}(0)=0.$
\end{description}

\noindent\textit{Moreover, in case (i)}
\[
\mathcal{L}_{\infty}(g|T)=\left[  \max_{i=1}^{K}\frac{\deg R_{i}}{i}\right]
^{-1}%
\]
\noindent\textit{and in case (iii)}
\[
\mathcal{L}_{\infty}(g|T)=-\left[  \min_{i=1}^{r}\frac{\operatorname*{ord}%
_{0}R_{i}}{r+1-i}\right]  ^{-1}.
\]

\begin{flushright}$\square$\end{flushright}

Let $f$ be a polynomial in two complex variables of the form (\ref{2.1}) and
$\deg f>1.$. Fix $\lambda_{0}\in\mathbb{C}$, denote $z:=(x,y)$ and define
\begin{align*}
S_{\lambda_{0}} &  :=\{z\in\mathbb{C}^{2}:f(z)=\lambda_{0}\},\\
S_{y} &  :=\{z\in\mathbb{C}^{2}:f_{y}^{\prime}(z)=0\}.
\end{align*}

In notation of Lemma \ref{L2.1} we put $\Phi_{i}(t):=(t^{D},\beta_{i}%
(t))$\ for $i\in\{1,...,n\}$ and as previously $\Psi_{j}(t):=(t^{D},\gamma
_{j}(t))$ for $j\in\{1,...,n-1\}.$

Under these notation we give, without a proof, a simple lemma which follows
easily from Lemma \ref{L2.1}.

\begin{lemma}
\label{L3.1}We have
\end{lemma}

(i) $\ \ \deg\Phi_{i}=D,\;i=1,...,n,\;\;\deg\Psi_{j}=D,\;j=1,...,n-1,$

(ii) $\ \mathcal{L}_{\infty}(f_{y}^{\prime}|S_{\lambda_{0}})=\frac{1}{D}%
\min_{i=1}^{n}\deg f_{y}^{\prime}\circ\Phi_{i},$

(iii) $\mathcal{L}_{\infty}(f-\lambda_{0}|S_{y})=\frac{1}{D}\min_{j=1}%
^{n-1}\deg\left(  f-\lambda_{0}\right)  \circ\Psi_{j}.$

\begin{flushright}$\square$\end{flushright}

Now, we give a theorem important in the sequel.

\begin{theorem}
\label{T3.2}If $\mathcal{L}_{\infty}(f-\lambda_{0},f_{y}^{\prime})<0$, then

\begin{description}
\item[(i)] $\mathcal{L}_{\infty}(f-\lambda_{0},f_{y}^{\prime})=\mathcal{L}%
_{\infty}(f-\lambda_{0}|S_{y}),$

\item[(ii)] $\lambda_{0}\in\Lambda(f).$
\end{description}

\noindent Moreover, if additionally $\mathcal{L}_{\infty}(f-\lambda_{0}%
,f_{y}^{\prime})\neq-\infty$ then
\begin{equation}
\mathcal{L}_{\infty}(f-\lambda_{0}|S_{y})<\mathcal{L}_{\infty}(f_{y}^{\prime
}|S_{\lambda_{0}}).\label{3.2}%
\end{equation}
\end{theorem}

\begin{proof}
Let us start from (i). In the case $\mathcal{L}_{\infty}(f-\lambda_{0}%
,f_{y}^{\prime})=-\infty$ we get easily (cf. \cite{CK4}, Theorem 3.1(iv)) that
$\mathcal{L}_{\infty}(f-\lambda_{0}|S_{y})=-\infty,$ which gives (i) in this case.

Let us assume now that $\mathcal{L}_{\infty}(f-\lambda_{0},f_{y}^{\prime}%
)\neq-\infty.$ In this case by the Main Theorem in \cite{CK1}, (cf.
\cite{CK3}, Theorem 1) we have
\begin{equation}
\mathcal{L}_{\infty}(f-\lambda_{0},f_{y}^{\prime})=\min(\mathcal{L}_{\infty
}(f-\lambda_{0}|S_{y}),\mathcal{L}_{\infty}(f_{y}^{\prime}|S_{\lambda_{0}%
})).\label{3.1}%
\end{equation}
Hence to prove (i) in this case it suffices to show (\ref{3.2})$.$

Assume to the contrary that (\ref{3.2}) does not hold. Then by (\ref{3.1}) and
the assumption of the theorem we have $\mathcal{L}_{\infty}(f_{y}^{\prime
}|S_{\lambda_{0}})<0.$ On the other hand, by Lemma \ref{L3.1}(ii) there exists
$i\in\{1,...,n\}$ such that
\begin{equation}
\deg f_{y}^{\prime}\circ\Phi_{i}=D\mathcal{L}_{\infty}(f_{y}^{\prime
}|S_{\lambda_{0}}).\label{3.3}%
\end{equation}
By the above we get $\deg f_{y}^{\prime}\circ\Phi_{i}<0.$ Hence and by Lemma
\ref{2.1}(c) we have
\begin{align*}
\deg f_{y}^{\prime}\circ\Phi_{i} &  =\sum_{j=1,j\neq i}^{n}\deg\left(
\beta_{i}-\beta_{j}\right) \\
&  >\sum_{j=1,j\neq i}^{n}\deg\left(  \beta_{i}-\beta_{j}\right)
+\min_{j=1,j\neq i}^{n}\deg(\beta_{i}-\beta_{j}).
\end{align*}
In consequence we get
\[
\deg f_{y}^{\prime}\circ\Phi_{i}>\min_{k=1}^{n}\left(  \sum_{j=1,j\neq k}%
^{n}\deg\left(  \beta_{k}-\beta_{j}\right)  +\min_{j=1,j\neq k}^{n}\deg
(\beta_{k}-\beta_{j})\right)  .
\]
Hence by Proposition \ref{P2.2}, Lemma \ref{L3.1}(iii) and (\ref{3.3})
\[
\mathcal{L}_{\infty}(f-\lambda_{0}|S_{y})<\mathcal{L}_{\infty}(f_{y}^{\prime
}|S_{\lambda_{0}}),
\]
which gives a contradiction. Then (\ref{3.2}) holds.

Assertion (ii) is a simple consequence of the facts $\mathcal{L}_{\infty
}(f-\lambda_{0}|S_{y})<0,$ Lemma \ref{LP} and (\ref{1.4}).

This ends the proof.
\end{proof}

Let us fix the same notation as in Theorem \ref{T3.2}.

\begin{theorem}
\label{T3.3}If $\mathcal{L}_{\infty}(f-\lambda_{0},f_{y}^{\prime})<0$, then
\begin{equation}
\mathcal{L}_{\infty,\lambda_{0}}(f)=\mathcal{L}_{\infty}(f-\lambda_{0}%
,f_{y}^{\prime})-1.\label{3.5}%
\end{equation}
\end{theorem}

\begin{proof}
Without loss of generality we may assume, as before, that $\mathcal{L}%
_{\infty}(f-\lambda_{0},f_{y}^{\prime})>-\infty.$ By Theorem \ref{T3.2}(i) it
suffices to show that
\begin{equation}
\mathcal{L}_{\infty,\lambda_{0}}(f)=\mathcal{L}_{\infty}(f-\lambda_{0}%
|S_{y})-1.\label{3.6}%
\end{equation}
We shall first show the inequality
\begin{equation}
\mathcal{L}_{\infty,\lambda_{0}}(f)\geqslant\mathcal{L}_{\infty}(f-\lambda
_{0}|S_{y})-1.\label{3.7}%
\end{equation}
According to definition (\ref{1.1}) of $\mathcal{L}_{\infty,\lambda_{0}}(f) $
it suffices to show that for any meromorphic curve $\Phi(t)=(\varphi
_{1}(t),\varphi_{2}(t))$ satisfying
\begin{gather}
\deg\Phi>0,\label{3.8}\\
\deg\left(  f-\lambda_{0}\right)  \circ\Phi<0,\label{3.9}%
\end{gather}
we have
\begin{equation}
\frac{\deg\nabla f\circ\Phi}{\deg\Phi}\geqslant\mathcal{L}_{\infty}%
(f-\lambda_{0}|S_{y})-1.\label{3.10}%
\end{equation}
From (\ref{3.8}) and (\ref{3.9}) it easily follows $\deg\varphi_{1}>0.$ So,
without loss of generality we may assume that $\Phi(t)=(t^{\deg\varphi_{1}%
},\varphi(t)).$ Then by (\ref{3.9}) we also get easily that $\deg\Phi
=\deg\varphi_{1}.$ On the other hand, by Lemma \ref{L3.1}(iii) it follows that
there exists $l_{\ast}\in\{1,...,n-1\}$ such that
\begin{equation}
\mathcal{L}_{\infty}(f-\lambda_{0}|S_{y})=\frac{\deg\left(  f-\lambda
_{0}\right)  \circ\Psi_{l_{\ast}}}{\deg\Psi_{l_{\ast}}}.\label{3.11}%
\end{equation}
Hence we get that inequality (\ref{3.10}) can be replaced by the inequality
\begin{equation}
\frac{\deg\nabla f\circ\Phi}{\deg\Phi}\geqslant\frac{\deg\left(  f-\lambda
_{0}\right)  \circ\Psi_{l_{\ast}}}{\deg\Psi_{l_{\ast}}}-1.\label{3.12}%
\end{equation}
At the cost of superpositions of $\Phi$ and $\Psi_{l_{\ast}}$, if necessary,
with appropriate powers of $t^{\alpha}$ and $t^{\beta},$ which does not change
the value of fraction in (\ref{3.12}), we may assume that $\deg\Phi=\deg
\Psi_{l_{\ast}}.$ Moreover, increasing $D$ in Lemma \ref{L2.1} we may also
assume that $\deg\Phi=D.$ Summing up, to show (\ref{3.7}) it suffices to
prove
\begin{equation}
\deg\nabla f\circ\Phi\geqslant\deg\left(  f-\lambda_{0}\right)  \circ
\Psi_{l_{\ast}}-D.\label{3.13}%
\end{equation}

Before the proof of this we notice that inequality (\ref{3.9}) implies easily
the following
\begin{equation}
\deg\left(  f-\lambda_{0}\right)  \circ\Phi\leqslant\deg\nabla f\circ
\Phi+D.\label{3.14.}%
\end{equation}

Consider now two cases:

(a) there exists $l_{0}\in\{1,...,n-1\}$ such that
\[
\deg(\varphi-\gamma_{l_{0}})<\min_{i=1}^{n}\deg(\varphi-\beta_{i}),
\]

(b) for each $l\in\{1,...,n-1\}$%
\[
\deg(\varphi-\gamma_{l})\geqslant\min_{i=1}^{n}\deg(\varphi-\beta_{i}).
\]
In the case (a) for each $j\in\{1,...,n\}$ we have
\[
\deg(\gamma_{l_{0}}-\beta_{j})=\deg(\gamma_{l_{0}}-\varphi+\varphi-\beta
_{j})=\deg(\varphi-\beta_{j}).
\]
Then
\[
\deg\left(  f-\lambda_{0}\right)  \circ\Psi_{l_{0}}\leqslant\deg\left(
f-\lambda_{0}\right)  \circ\Phi.
\]
Hence, from (\ref{3.11}) and Lemma \ref{L3.1}(iii) we get
\begin{equation}
\deg\left(  f-\lambda_{0}\right)  \circ\Psi_{l_{\ast}}\leqslant\deg\left(
f-\lambda_{0}\right)  \circ\Phi.\label{3.15}%
\end{equation}
By (\ref{3.14.}) and (\ref{3.15}) we get (\ref{3.13}) in case (a).

We shall now show (\ref{3.13}) in the case (b). Let $\min_{i=1}^{n}%
\deg(\varphi-\beta_{i})=\deg(\varphi-\beta_{i_{0}})$ for some $i_{0}%
\in\{1,...,n\}.$ Then for each $l\in\{1,...,n-1\}$%
\[
\deg(\beta_{i_{0}}-\gamma_{l})=\deg(\beta_{i_{0}}-\varphi+\varphi-\gamma
_{l})\leqslant\deg(\varphi-\gamma_{l}).
\]
Hence
\begin{equation}
\deg f_{y}^{\prime}\circ\Phi_{i_{0}}\leqslant\deg f_{y}^{\prime}\circ
\Phi.\label{3.16}%
\end{equation}
On the other hand, by Proposition \ref{P2.2}, Lemma \ref{L3.1}(iii), and
(\ref{3.11})
\begin{align*}
\deg f_{y}^{\prime}\circ\Phi_{i_{0}} &  =\sum_{j=1,j\neq i_{0}}^{n}\deg\left(
\beta_{i_{0}}-\beta_{j}\right) \\
&  =\sum_{j=1,j\neq i_{0}}^{n}\deg\left(  \beta_{i_{0}}-\beta_{j}\right)
+\min_{j=1,j\neq i_{0}}^{n}\deg(\beta_{i_{0}}-\beta_{j})-\min_{j=1,j\neq
i_{0}}^{n}\deg(\beta_{i_{0}}-\beta_{j})\\
&  \geqslant\min_{k=1}^{n}\left(  \sum_{j=1,j\neq k}^{n}\deg\left(  \beta
_{k}-\beta_{j}\right)  +\min_{j=1,j\neq k}^{n}\deg(\beta_{k}-\beta
_{j})\right)  -D\\
&  =D\mathcal{L}_{\infty}(f-\lambda_{0}|S_{y})-D=\deg\left(  f-\lambda
_{0}\right)  \circ\Psi_{l_{1}}-D.
\end{align*}
Hence and by (\ref{3.16}) we get
\begin{equation}
\deg f_{y}^{\prime}\circ\Phi\geqslant\deg\left(  f-\lambda_{0}\right)
\circ\Psi_{l_{1_{\ast}}}-D.\label{3.17}%
\end{equation}
By (\ref{3.17}) and the obvious inequality $\deg\nabla f\circ\Phi\geqslant\deg
f_{y}^{\prime}\circ\Phi$ we get inequality (\ref{3.13}) in case (b). Then we
have proved (\ref{3.13}) and in consequence (\ref{3.7}).

To finish the proof we have to show
\begin{equation}
\mathcal{L}_{\infty,\lambda_{0}}(f)\leqslant\mathcal{L}_{\infty}(f-\lambda
_{0}|S_{y})-1.\label{3.18}%
\end{equation}
By assumption, Theorem \ref{T3.2}(i) and (\ref{3.11}) we have
\begin{equation}
\deg\left(  f-\lambda_{0}\right)  \circ\Psi_{l_{\ast}}<0.\label{3.19}%
\end{equation}
Hence
\begin{equation}
\deg\left(  f-\lambda_{0}\right)  \circ\Psi_{l_{\ast}}=\deg\nabla f\circ
\Psi_{l_{\ast}}+D.\label{3.20}%
\end{equation}
Hence and from (\ref{3.11}) we get
\[
\mathcal{L}_{\infty}(f-\lambda_{0}|S_{y})-1=\frac{\deg\nabla f\circ
\Psi_{l_{\ast}}}{\deg\Psi_{l_{\ast}}}.
\]
Hence taking into account (\ref{3.19}) and (\ref{1.1}) we obtain (\ref{3.18}).

This ends the proof of the theorem.
\end{proof}

We shall now give three simple corollaries of Theorems \ref{T3.2} and
\ref{T3.3}.

\begin{corollary}
\label{C3.4}The following conditions are equivalent:
\end{corollary}

(i) $\ \ \mathcal{L}_{\infty}(f-\lambda_{0},f_{y}^{\prime})<0,$

(ii) $\ \mathcal{L}_{\infty,\lambda_{0}}(f)<-1,$

(iii) $\mathcal{L}_{\infty,\lambda_{0}}(f)<0,$

(iv) $\lambda_{0}\in\Lambda(f),$

(v) \ $\mathcal{L}_{\infty}(f-\lambda_{0}|S_{y})<0.$

\begin{proof}
(i)$\Rightarrow$(ii). Theorem \ref{T3.3}.

(ii)$\Rightarrow$(iii). Obvious.

(iii)$\Rightarrow$(i). Obvious.

(i)$\Rightarrow$(iv). Theorem \ref{T3.2}(ii).

(iv)$\Rightarrow$(i). \cite{CK4}, Theorem 3.1.

(i)$\Rightarrow$(v). Theorem \ref{T3.2}(i).

(v)$\Rightarrow$(i). Obvious.

This ends the proof.
\end{proof}

\begin{corollary}
\label{C3.5}If $\mathcal{L}_{\infty}(\nabla f)\leqslant-1,$ then

\begin{description}
\item[(i)] there exists $\lambda_{0}\in\mathbb{C}$ such that $\mathcal{L}%
_{\infty}(\nabla f)=\mathcal{L}_{\infty,\lambda_{0}}(f),$

\item[(ii)] $\mathcal{L}_{\infty}(\nabla f)=\mathcal{L}_{\infty}(\nabla
f|S_{y}).$
\end{description}
\end{corollary}

\begin{proof}
Let $\Phi,$ $\deg\Phi>0,$ be a meromorphic curve on which the \L ojasiewicz
exponent $\mathcal{L}_{\infty}(\nabla f)$ is attained. Then
\begin{equation}
\mathcal{L}_{\infty}(\nabla f)=\frac{\deg\nabla f\circ\Phi}{\deg\Phi
}.\label{3.21}%
\end{equation}
We shall show
\begin{equation}
\deg f\circ\Phi\leqslant0.\label{3.22}%
\end{equation}
Indeed, it suffices to consider the case $\deg f\circ\Phi\neq0.$ Then
\[
\frac{\deg f\circ\Phi}{\deg\Phi}\leqslant\frac{\deg\nabla f\circ\Phi}{\deg
\Phi}+1=\mathcal{L}_{\infty}(\nabla f)+1\leqslant0,
\]
which gives (\ref{3.22}).

Inequality (\ref{3.22}) implies that there exists $\lambda_{0}\in\mathbb{C} $
such that
\begin{equation}
\deg(f-\lambda_{0})\circ\Phi<0.\label{3.23}%
\end{equation}
Then by (\ref{3.21}), (\ref{3.23}) and (\ref{1.1}) we get $\mathcal{L}%
_{\infty,\lambda_{0}}(f)\leqslant\mathcal{L}_{\infty}(\nabla f).$ The opposite
inequality is obvious. This gives (i).

From (\ref{3.23}), the assumption and (\ref{3.21}) we get $\mathcal{L}%
_{\infty}(f-\lambda_{0},f_{y}^{\prime})<0.$ Hence according to (i) and
Theorems \ref{T3.3} and \ref{T3.2}(i) we get
\[
\mathcal{L}_{\infty}(\nabla f)=\mathcal{L}_{\infty}(f-\lambda_{0}|S_{y})-1.
\]
Hence and from the obvious inequality
\[
\mathcal{L}_{\infty}(f-\lambda_{0}|S_{y})-1\geqslant\mathcal{L}_{\infty
}(\nabla f|S_{y})
\]
we obtain
\[
\mathcal{L}_{\infty}(\nabla f)\geqslant\mathcal{L}_{\infty}(\nabla f|S_{y}).
\]
The opposite inequality is obvious, which gives (ii).

This ends the proof.
\end{proof}

\begin{corollary}
\label{C3.6}If $\mathcal{L}_{\infty}(f-\lambda_{0},f_{y}^{\prime})<0$ and
functions $\beta_{1},...,\beta_{n},$ meromorphic at infinity, are as in Lemma
\ref{L2.1} then
\begin{equation}
\mathcal{L}_{\infty,\lambda_{0}}(f)+1=\frac{1}{D}\min_{i=1}^{n}\left(
\sum_{j=1,j\neq i}^{n}\deg\left(  \beta_{i}-\beta_{j}\right)  +\min_{j=1,j\neq
i}^{n}\deg(\beta_{i}-\beta_{j})\right)  .\label{3.24}%
\end{equation}
\end{corollary}

\begin{proof}
By Theorems \ref{T3.3} and \ref{T3.2} (i) we get
\begin{equation}
\mathcal{L}_{\infty,\lambda_{0}}(f)+1=\mathcal{L}_{\infty}(f-\lambda_{0}%
|S_{y}).\label{3.25}%
\end{equation}
According to Lemma \ref{L3.1}
\begin{equation}
\mathcal{L}_{\infty}(f-\lambda_{0}|S_{y})=\frac{1}{D}\min_{l=1}^{n-1}%
\deg(f-\lambda_{0})\circ\Psi_{l}.\label{3.26}%
\end{equation}
Now, (\ref{3.25}), (\ref{3.26}) and Proposition \ref{P2.2} implies (\ref{3.24}).

This ends the proof.
\end{proof}

At the end of this section we notice that from Corollary \ref{C3.4} it follows
that all results of this section concern critical values of $f$ at infinity.
Indeed, by Corollary \ref{C3.4} one can always replace the assumption
$\mathcal{L}_{\infty}(f-\lambda_{0},f_{y}^{\prime})<0$ with the assumption
$\lambda_{0}\in\Lambda(f).$

We shall now discuss the connection of the above three corollaries with the
results by Ha \cite{HA2}. It shall be shown in Section 5 that the above \L
ojasiewicz exponent $\mathcal{L}_{\infty,\lambda}(f),$ defined by (\ref{1.1}),
coincides with the \L ojasiewicz exponent $\widetilde{\mathcal{L}}%
_{\infty,\lambda}(f)$, defined by (\ref{1.1'}), introduced by Ha in
\cite{HA2}. Thus Corollary \ref{C3.4} is a changed and extended version of
Theorems 1.3.1 and 1.3.2 in \cite{HA2}. In turn, Corollaries \ref{C3.5}(i) and
\ref{C3.6} correspond exactly to Theorems 1.4.3 and 1.4.1 in \cite{HA2}, respectively.

\section{Effective formulas for $\mathcal{L}_{\infty,\lambda}(f)$}

In this section $f$ is a polynomial in two complex variables of the form
(\ref{2.1}). Let $(\lambda,u)\in\mathbb{C}^{2}$ and $Q(x,\lambda
,u):=\operatorname*{Res}_{y}(f-\lambda,f_{y}^{\prime}-u)$ be the resultant of
the polynomials $f-\lambda$ and $f_{y}^{\prime}-u$ with respect to the
variable $y.$ By the definition of the resultant we get easily that
$Q(0,\lambda,0)=\pm n^{n}\lambda^{n-1}+$ terms of lower degrees. Hence
$Q\neq0.$ We put
\begin{equation}
Q(x,\lambda,u)=Q_{0}(\lambda,u)x^{N}+\cdots+Q_{N}(\lambda,u),\;\;\;Q_{0}%
\neq0.\label{4.1}%
\end{equation}

Let us pass now to the effective calculations of $\mathcal{L}_{\infty,\lambda
}(f).$ We start with the first main theorem conerning the case when
$\lambda_{0}$ is a critical value of $f$ at infinity.

\begin{theorem}
\label{T4.1}A point $\lambda_{0}\in\mathbb{\mathbb{C}}$ is a critical value of
$f$ at infinity if and only if $Q_{0}(\lambda_{0},0)=0.$ Moreover

\begin{description}
\item[(i)] $\mathcal{L}_{\infty,\lambda_{0}}(f)=-\infty$ if and only if
$Q_{0}(\lambda_{0},0)=\cdots=Q_{N}(\lambda_{0},0)=0,$

\item[(ii)]
\[
\mathcal{L}_{\infty,\lambda_{0}}(f)=-1-\left[  \min_{i=0}^{r}%
\frac{\operatorname*{ord}_{(\lambda_{0},0)}Q_{i}(\lambda,u)}{r+1-i}\right]
^{-1}%
\]
if and only if there exists $r\in\{0,...,N-1\}$ such that $Q_{0}(\lambda
_{0},0)=\cdots=Q_{r}(\lambda_{0},0)=0,$ $Q_{r+1}(\lambda_{0},0)\neq0.$
\end{description}
\end{theorem}

\begin{proof}
By Corollary \ref{C3.4} (iv)$\iff$(i) and Theorem 3.1 in \cite{CK4} we get the
first assertion of the theorem. The second one follows from Theorems 3.1 and
3.3 in \cite{CK4} which give \ formulas for $\mathcal{L}_{\infty}%
(f-\lambda_{0},f_{y}^{\prime})$ in the case $Q_{0}(\lambda_{0},0)=0$ and
Theorem \ref{T3.3}.
\end{proof}

Now let us pass to the more complicated case when $\lambda_{0}$ is not a
critical value of $f$ at infinity$.$ Similarly as previously we shall use
resultant (\ref{4.1}). We shall now prove two theorems on dependence of
$\mathcal{L}_{\infty,\lambda_{0}}(f)$ for $\lambda_{0}\notin\Lambda(f)$ on the
coefficient $Q_{0}(\lambda,u)$ of the resultant (\ref{4.1}).

Directly from the first assertion of the above theorem we have.

\begin{corollary}
\label{L4.3}A point $\lambda_{0}\notin\Lambda(f)$ if and only if
$Q_{0}(\lambda_{0},0)\neq0.$
\end{corollary}

\begin{flushright}$\square$\end{flushright}

The first theorem shall be preceded by a lemma, well known in local case.
First we introduce notations.

Let $\mathcal{M}(t)$ be the field of germs of meromorphic functions at
infinity i.e. the field of all Laurent series of the form $\sum_{n=k}%
^{-\infty}a_{n}t^{n},$ $k\in\mathbb{Z},$ convergent in a neighbourhood of
$\infty\in\overline{\mathbb{\mathbb{C}}}$. Let $\mathcal{M}(t)^{\ast}%
:=\bigcup_{k=1}^{\infty}\mathcal{M}(t^{1/k})$ be the field of convergent
Puiseux series at infinity. Similarly as in the local case $\mathcal{M}%
(t)^{\ast}$ is an algebraically closed field. If $\varphi\in\mathcal{M}%
(t)^{\ast}$ and $\varphi(t)=\psi(t^{1/k})$ for $\psi\in\mathcal{M}(t)$, then
we define $\deg\varphi:=(1/k)\deg\psi.$

Using simple properties of the function $\deg$ and the Vieta formulae we obtain

\begin{lemma}
\label{LA.1}Let
\[
P(x,t)=c_{0}(t)x^{N}+c_{1}(t)x^{N-1}+\cdots+c_{N}(t)=c_{0}(t)(x-\varphi
_{1}(t))\cdots(x-\varphi_{N}(t)),
\]
where $c_{0},c_{1},...,c_{N}\in\mathcal{M}(t),$ $c_{0}\neq0,$ $\varphi
_{1},...,\varphi_{N}\in\mathcal{M}(t)^{\ast}.$ Then
\[
\max_{i=1}^{N}\deg\varphi_{i}=\max_{i=1}^{N}\frac{\deg c_{i}-\deg c_{0}}{i}.
\]
\end{lemma}

We shall now prove the second main theorem of the paper.

\begin{theorem}
\label{T4.7}If $\lambda_{0}\notin\Lambda(f)$ and $\deg_{u}Q_{0}(\lambda,u)=0$
then $\mathcal{L}_{\infty,\lambda_{0}}(f)>0$ and
\begin{equation}
\mathcal{L}_{\infty,\lambda_{0}}(f)=\left[  \max_{i=1}^{N}\frac{\deg_{u}%
Q_{i}(\lambda,u)}{i}\right]  ^{-1}.\label{4.6}%
\end{equation}
\end{theorem}

\begin{proof}
Put $\delta:=\left[  \max_{i=1}^{N}\frac{\deg_{u}Q_{i}(\lambda,u)}{i}\right]
^{-1}$. By an elementary property of the resultant $Q$ it follows that
$\deg_{u}Q(\lambda,u)>0.$ Hence $\delta>0.$

We first show that
\begin{equation}
\delta\leqslant\mathcal{L}_{\infty,\lambda_{0}}(f).\label{4.7}%
\end{equation}
Take an arbitrary meromorphic curve $\Phi(t)=(x(t),y(t))$ such that $\deg
\Phi>0$ and $\deg(f-\lambda_{0})\circ\Phi<0.$ To show (\ref{4.7}) it suffices
to prove
\begin{equation}
\frac{\deg f_{y}^{\prime}\circ\Phi}{\deg\Phi}\geqslant\delta.\label{4.8}%
\end{equation}
Notice that inequality $\deg(f-\lambda_{0})\circ\Phi<0$ and (\ref{2.1})
implies immediately $\deg\Phi=\deg x.$ Put $\lambda(t):=f\circ\Phi(t),$
$u(t):=f_{y}^{\prime}\circ\Phi(t).$ Then (\ref{4.8}) takes the form
\begin{equation}
\frac{\deg u}{\deg x}\geqslant\delta.\label{4.9}%
\end{equation}

By a property of the resultant we have
\begin{equation}
Q(x(t),\lambda(t),u(t))\equiv0.\label{4.10}%
\end{equation}
By Corollary \ref{L4.3} and the assumption of the theorem we have
$Q_{0}(\lambda_{0},0)\neq0$ and $Q_{0}$ does not depend on $u.$ Since
$\deg(\lambda(t)-\lambda_{0})<0$ then
\begin{equation}
\deg Q_{0}(\lambda(t),u(t))=0.\label{4.11}%
\end{equation}
By (\ref{4.11}) and (\ref{4.10}) taking into account $\deg x>0$ and
$\deg\lambda\leqslant0$ we get easily
\begin{equation}
\deg u>0.\label{4.12}%
\end{equation}
Consider the polynomial in variable $x$
\[
Q(x,\lambda(t),u(t))=Q_{0}(\lambda(t),u(t))x^{N}+\cdots+Q_{N}(\lambda(t),u(t))
\]
with coefficients meromorphic at infinity. Identifying meromorphic functions
at infinity with with their germs in $\mathcal{M}(t)$ and using (\ref{4.10}),
(\ref{4.11}) and Lemma \ref{LA.1} we get
\[
\deg x(t)\leqslant\max_{i=1}^{N}\frac{\deg Q_{i}(\lambda(t),u(t))}{i}.
\]
Hence and from the inequality $\deg\lambda(t)\leqslant0$ we obtain
\[
\deg x(t)\leqslant\deg u(t)\max_{i=1}^{N}\frac{\deg_{u}Q_{i}(\lambda,u)}%
{i}=\frac{1}{\delta}\deg u(t).
\]

This gives (\ref{4.9}) and in consequence (\ref{4.8}) and then (\ref{4.7}).

We shall now prove that
\begin{equation}
\mathcal{L}_{\infty,\lambda_{0}}(f)\leqslant\delta.\label{4.16}%
\end{equation}
Let us introduce notations. Let
\begin{align}
Q_{i}(\lambda,u) &  =Q_{0}^{i}(\lambda)u^{k_{i}}+\cdots+Q_{k_{i}}^{i}%
(\lambda),\label{4.17}\\
\alpha_{i} &  :=\operatorname*{ord}\nolimits_{\lambda_{0}}Q_{0}^{i}%
(\lambda)\label{4.18}%
\end{align}
for $i=1,...,N.$ Take now an arbitrary $M\in\mathbb{N}$ and put
\begin{equation}
u_{M}(t):=t^{M},\;\;\;\;\;\lambda(t):=\lambda_{0}+\frac{1}{t}.\label{4.19}%
\end{equation}
Since $Q_{0}(\lambda_{0},0)\neq0$ and $Q_{0}$ does not depend on $u$ then
\begin{equation}
\deg Q_{0}(\lambda(t),u_{M}(t))=0.\label{4.11M}%
\end{equation}
On the other hand by (\ref{4.17}), (\ref{4.18}) and (\ref{4.19}) for almost
all $M\in\mathbb{N}$ we have
\begin{equation}
\deg Q_{i}(\lambda(t),u_{M}(t))=Mk_{i}-\alpha_{i},\;\;\;i=1,...,N.\label{A.16}%
\end{equation}
Define $\alpha:=\min\{\alpha_{i}:\frac{k_{i}}{i}=\delta^{-1}\}.$ Then from
(\ref{A.16}) and definition of $\delta$ it follows that for almost all
$M\in\mathbb{N}$ we have
\begin{equation}
\max_{i=1}^{N}\frac{\deg Q_{i}(\lambda(t),u_{M}(t))}{i}=M\delta^{-1}%
-\alpha.\label{A.17}%
\end{equation}

Consider the polynomial in variable $x$
\[
Q(x,\lambda(t),u_{M}(t))=Q_{0}(\lambda(t),u_{M}(t))x^{N}+\cdots+Q_{N}%
(\lambda(t),u_{M}(t))
\]
with coefficients meromorphic at infinity. Identifying meromorphic functions
at infinity with their germs in $\mathcal{M}(t),$ applying Lemma \ref{LA.1}
and taking into account (\ref{4.11M}), we obtain that there exists
$\varphi_{M}\in\mathcal{M}(t)^{\ast}$ such that
\begin{equation}
Q(\varphi_{M}(t),\lambda(t),u_{M}(t))=0\;\;\;\text{in }\mathcal{M}(t)^{\ast
}\label{A.19}%
\end{equation}
and
\[
\deg\varphi_{M}(t)=\max_{i=1}^{N}\frac{\deg Q_{i}(\lambda(t),u_{M}(t))}{i}.
\]

Taking into account (\ref{A.17}) for almost all $M\in\mathbb{N}$ we have
\begin{equation}
\deg\varphi_{M}(t)=M\delta^{-1}-\alpha.\label{A.20}%
\end{equation}

Consider polynomials $f(\varphi_{M}(t),y)-\lambda(t)$ and $f_{y}^{\prime
}(\varphi_{M}(t),y)-u_{M}(t)$ in variable $y$ with coefficients in
$\mathcal{M}(t)^{\ast}.$ By (\ref{A.19}) their resultant is equal to zero in
$\mathcal{M}(t)^{\ast}.$ Then there exists a function $\psi_{M}(t)\in
\mathcal{M}(t)^{\ast}$ such that
\begin{align}
f(\varphi_{M}(t),\psi_{M}(t)) &  =\lambda(t),\label{A.21}\\
f_{y}^{\prime}(\varphi_{M}(t),\psi_{M}(t)) &  =u_{M}(t).\label{A.22}%
\end{align}

By (\ref{A.21}) and (\ref{2.1}) we get
\begin{equation}
\deg\varphi_{M}\geqslant\deg\psi_{M}.\label{A.23}%
\end{equation}

By definition of $\mathcal{M}(t)^{\ast}$ there exists $D\in\mathbb{N}$ such
that $\varphi_{M}(t^{D}),\psi_{M}(t^{D})\in\mathcal{M}(t).$ Then there exists
a meromorphic mapping at infinity $\Phi_{M}(t)=(x_{M}(t),y_{M}(t))$ such that
the germs of $x_{M}(t)$ and $y_{M}(t)$ at infinity are equal $\varphi
_{M}(t^{D})$ and $\psi_{M}(t^{D}),$ respectively. By (\ref{A.23}) and
(\ref{A.20}) we obtain
\begin{equation}
\deg\Phi_{M}=\deg x_{M}=D(M\delta^{-1}-\alpha).\label{A.24}%
\end{equation}
By (\ref{A.21}) and (\ref{A.22}) we have
\begin{align}
\deg(f-\lambda_{0})\circ\Phi_{M} &  =\deg(\lambda(t^{D})-\lambda
_{0})=-D<0,\label{A.25}\\
\deg f_{y}^{\prime}\circ\Phi_{M} &  =\deg u_{M}(t^{D})=MD.\label{A.26}%
\end{align}

By (\ref{A.24}) for almost all $M\in\mathbb{N}$ $\deg\Phi_{M}>0.$ Hence, from
(\ref{A.25}) and the first equality in (\ref{A.24}) we easily get $\deg
f_{x}^{\prime}\circ\Phi_{M}\leqslant\deg f_{y}^{\prime}\circ\Phi_{M},$ that
is
\begin{equation}
\deg\nabla f\circ\Phi_{M}=\deg f_{y}^{\prime}\circ\Phi_{M}.\label{A.27}%
\end{equation}

By (\ref{A.24}), (\ref{A.26}) and (\ref{A.27}) for almost all $M\in\mathbb{N}$
we get
\[
\frac{\deg\nabla f\circ\Phi_{M}}{\deg\Phi_{M}}=\frac{M}{M\delta^{-1}-\alpha}.
\]

Summing up, we have found a sequence $\{\Phi_{M}\}$ of mappings meromorphic at
infinity, such that for almost all $M\in\mathbb{N}$, $\deg\Phi_{M}>0,$
(\ref{A.25}) holds and
\[
\lim_{M\rightarrow\infty}\frac{\deg\nabla f\circ\Phi_{M}}{\deg\Phi_{M}}%
=\delta,
\]
which in according with (\ref{1.1}) gives (\ref{4.16}).

This ends the proof.
\end{proof}

We shall now prove the third main theorem of the paper.

\begin{theorem}
\label{T4.8}If $\lambda_{0}\notin\Lambda(f)$ and $\deg_{u}Q_{0}(\lambda,u)>0$
then $\mathcal{L}_{\infty,\lambda_{0}}(f)=0.$
\end{theorem}

\begin{proof}
If $\deg f=1$ then we easily check that $\mathcal{L}_{\infty,\lambda_{0}%
}(f)=0.$ Assume that $\deg f>1.$ Let us notice first that by assumption and
Corollary \ref{C3.4} (iii)$\Rightarrow$(iv)
\begin{equation}
\mathcal{L}_{\infty,\lambda_{0}}(f)\geqslant0.\label{4.30}%
\end{equation}

On the other hand, by assumption and Corollary \ref{L4.3} we have
\begin{equation}
Q_{0}(\lambda_{0},0)\neq0.\label{4.31}%
\end{equation}
To show the inequality opposite to (\ref{4.30}) we consider two cases:

(a) $\deg Q_{0}(\lambda_{0},u)>0,$

(b) $\deg Q_{0}(\lambda_{0},u)=0.$

In case (a) taking into account (\ref{4.31}) and Lemma \ref{LP} we get
$\mathcal{L}_{\infty}(f_{y}^{\prime}|f^{-1}(\lambda_{0}))=0.$ Hence and by
Lemma \ref{L3.1} (ii) it follows that there exists $i\in\{1,...,n\}$ such that
$\deg f_{y}^{\prime}\circ\Phi_{i}=0.$ On the other hand, from the form of $f$
we easily conclude that $\deg\nabla f\circ\Phi_{i}=\deg f_{y}^{\prime}%
\circ\Phi_{i}.$ Summing up, $\deg\Phi_{i}>0$ and
\begin{equation}
\deg(f-\lambda_{0})\circ\Phi_{i}=-\infty,\;\;\;\;\;\deg\nabla f\circ\Phi
_{i}=0.\label{4.32}%
\end{equation}
Then definition (\ref{1.1}) of $\mathcal{L}_{\infty,\lambda_{0}}(f)$ implies
$\mathcal{L}_{\infty,\lambda_{0}}(f)\leqslant0,$ which together with
(\ref{4.30}) gives the assertion of the theorem in case (a).

Let us pass to case (b).

We first show that there exist functions $\lambda(t),u(t),$ meromorphic at
infinity, such that
\begin{align}
&  \deg(\lambda(t)-\lambda_{0})<0,\label{B.1}\\
&  \deg u(t)>0,\label{B.2}\\
&  Q_{0}(\lambda(t),u(t))\equiv0.\label{B.3}%
\end{align}
Let
\[
Q_{0}(\lambda,u)=Q_{0}^{0}(\lambda)u^{k_{0}}+\cdots+Q_{k_{0}}^{0}%
(\lambda),\;\;Q_{0}^{0}\neq0.
\]

\noindent In this case taking into account (\ref{4.31}) we have
\begin{equation}
\operatorname*{ord}\nolimits_{\lambda_{0}}Q_{k_{0}}^{0}%
=0,\;\;\;\;\;\operatorname*{ord}\nolimits_{\lambda_{0}}Q_{i}^{0}%
>0\;\;\;\text{for \ }i=0,...,k_{0}-1.\label{B.4}%
\end{equation}
Put $\tilde{\lambda}(t):=\lambda_{0}+\frac{1}{t}$ and consider the polynomial
in variable $u$%
\begin{equation}
Q_{0}(\tilde{\lambda}(t),u)=Q_{0}^{0}(\tilde{\lambda}(t))u^{k_{0}}%
+\cdots+Q_{k_{0}}^{0}(\tilde{\lambda}(t))\label{B.5}%
\end{equation}
with coefficients meromorphic at infinity. Notice first that by (\ref{B.4})
and definition of $\tilde{\lambda}$%
\begin{align}
&  \max_{i=1}^{k_{0}}\frac{\deg Q_{i}^{0}(\tilde{\lambda}(t))-\deg Q_{0}%
^{0}(\tilde{\lambda}(t))}{i}\geqslant\label{B.6}\\
&  \geqslant\frac{\deg Q_{k_{0}}^{0}(\tilde{\lambda}(t))-\deg Q_{0}^{0}%
(\tilde{\lambda}(t))}{k_{0}}=\frac{1}{k_{0}}\operatorname*{ord}%
\nolimits_{\lambda_{0}}Q_{0}^{0}>0.\nonumber
\end{align}
Identifying meromorphic functions at infinity with their germs in
$\mathcal{M}(t)$ and using Lemma \ref{LA.1} to polynomial (\ref{B.5}) and
taking into account (\ref{B.6}) we get that there exists $\varphi
\in\mathcal{M}(t)^{\ast}$ such that
\begin{gather}
Q_{0}(\tilde{\lambda}(t),\varphi(t))=0\;\;\;\text{in }\mathcal{M}(t)^{\ast
},\label{B.7}\\
\deg\varphi>0.\label{B.8}%
\end{gather}
By definition of $\mathcal{M}(t)^{\ast}$ there exists $D\in\mathbb{N}$ such
that $\varphi(t^{D})\in\mathcal{M}(t).$ Then there exists a function $u(t)$
meromorphic at infinity which germs at infinity is equal to $\varphi(t^{D}).$
Put $\lambda(t):=\tilde{\lambda}(t^{D}).$ By definition of $\tilde{\lambda},$
(\ref{B.8}) and (\ref{B.7}) we obtain that the functions $\lambda(t),u(t)$
satisfy (\ref{B.1}), (\ref{B.2}) and (\ref{B.3}).

Now we consider two subcases:

(b$_{1}$) for every $i\in\{1,...,N\}$ \ $Q_{i}(\lambda(t),u(t))\equiv0,$

(b$_{2}$) there exists $i\in\{1,...,N\}$ such that \ $Q_{i}(\lambda
(t),u(t))\not \equiv 0.$

Consider case (b$_{1}$). Then for every $M\in\mathbb{N}$ we have
\begin{equation}
Q(t^{M},\lambda(t),u(t))\equiv0.\label{B.9}%
\end{equation}
Consider polynomials in variable $y$,
\[
f(t^{M},y)-\lambda(t),\;\;\;\;\;f_{y}^{\prime}(t^{M},y)-u(t)
\]
with coefficients meromorphic at infinity. According to (\ref{B.9}) their
resultant vanishes identically. Then there exist a function $y_{M}(t),$
meromorphic at infinity, and an integer $D\in\mathbb{N}$ such that
\begin{align}
f(t^{MD},y_{M}(t)) &  =\lambda(t^{D}),\label{B.10}\\
f_{y}^{\prime}(t^{MD},y_{M}(t)) &  =u(t^{D}).\label{B.11}%
\end{align}
Define $\Phi_{M}(t):=(t^{MD},y_{M}(t)).$ By (\ref{B.10}) and (\ref{B.1}) we
get
\begin{equation}
\deg(f\circ\Phi_{M}-\lambda_{0})<0.\label{B.12}%
\end{equation}
Hence and from (\ref{2.1}) we obtain
\begin{equation}
\deg\Phi_{M}=MD.\label{B.13}%
\end{equation}
By (\ref{B.12}) and (\ref{B.13}) we easily get $\deg f_{x}^{\prime}\circ
\Phi_{M}\leqslant\deg f_{y}^{\prime}\circ\Phi_{M},$ which means $\deg\nabla
f\circ\Phi_{M}=\deg f_{y}^{\prime}\circ\Phi_{M}.$ Hence by (\ref{B.11}) we
obtain
\begin{equation}
\frac{\deg\nabla f\circ\Phi_{M}}{\deg\Phi_{M}}=\frac{\deg u}{M}.\label{B.14}%
\end{equation}

Summing up, we have found a sequence $\{\Phi_{M}\}$ of mappings meromorphic at
infinity, such that $\deg\Phi_{M}>0$ and (\ref{B.12}) holds. Moreover, by
(\ref{B.14})
\[
\lim_{M\rightarrow\infty}\frac{\deg\nabla f\circ\Phi_{M}}{\deg\Phi_{M}}=0,
\]
which according to definition (\ref{1.1}) gives $\mathcal{L}_{\infty
,\lambda_{0}}(f)\leqslant0.$ Hence and by (\ref{4.30}) we get the assertion of
the theorem in case (b$_{1}$).

Let us pass to case (b$_{2}$). Put $\lambda_{M}(t):=\lambda(t)+\frac{1}{t^{M}%
}$ for $M\in\mathbb{N}$ and $M>-\deg(\lambda(t)-\lambda_{0}).$ It is easy to
see that for every $i\in\{0,1,...,N\}$ there exists an integer $\alpha_{i}%
\in\mathbb{Z}$ such that for any $M$%
\begin{equation}
\deg(Q_{i}(\lambda(t),u(t))-Q_{i}(\lambda_{M}(t),u(t)))\leqslant-M+\alpha
_{i}.\label{B.15}%
\end{equation}
Indeed, let $Q_{i}(\lambda,u(t))=R_{0}^{i}(t)\lambda^{l_{i}}+\cdots+R_{l_{i}%
}^{i}(t).$ Then it suffices to take $\alpha_{i}:=\max_{j=0}^{l_{i}}\deg
R_{j}^{i}(t).$ By (\ref{B.15}) we also easily obtain that for almost all $M$
there exists an integer $d_{M}\in\mathbb{Z}$ such that
\begin{gather}
\deg Q_{0}(\lambda_{M}(t),u(t)))=d_{M},\label{B.16}\\
d_{M}\leqslant-M+\alpha_{0}.\label{B.17}%
\end{gather}
Indeed, in our case (b) by (\ref{4.31}) and (\ref{B.3}) we see that the
polynomial $Q_{0}(\lambda,u(t))$ in variable $\lambda$ is not constant i.e.
$\deg_{\lambda}Q_{0}(\lambda,u(t))>0$. Then for almost all $M$ there exists an
integer $d_{M}\in\mathbb{Z}$ such that (\ref{B.16}) holds. Inequality
(\ref{B.17}) is a conequence of (\ref{B.15}) and (\ref{B.16}).

Consider now for sufficiently large $M$ the following polynomial in variable
$x$
\begin{equation}
Q(x,\lambda_{M}(t),u(t)))=Q_{0}(\lambda_{M}(t),u(t))x^{N}+\cdots+Q_{N}%
(\lambda_{M}(t),u(t)),\label{B.18}%
\end{equation}
with coefficients meromorphic at infinity.

Let us notice there exists $i_{0}\in\{1,...,N\}$ such that
\begin{equation}
\lim_{M\rightarrow\infty}(\deg Q_{i_{0}}(\lambda_{M}(t),u(t))-\deg
Q_{0}(\lambda_{M}(t),u(t)))=\infty.\label{B.19}%
\end{equation}
In fact, by assumption (b$_{2}$) there exists $i_{0}$ such that $Q_{i_{0}%
}(\lambda(t),u(t)))\not \equiv 0.$ Put $d_{i_{0}}:=\deg Q_{i_{0}}%
(\lambda(t),u(t))).$ Then $d_{i_{0}}\in\mathbb{Z}$ and by (\ref{B.15}) for
$-M<d_{i_{0}}-\alpha_{i_{0}}$ we have $\deg Q_{i_{0}}(\lambda_{M}%
(t),u(t))=d_{i_{0}}.$ Hence by (\ref{B.16}) and (\ref{B.17}) we get
(\ref{B.19}).

Let us return to polynomial (\ref{B.18}). Identifying meromorphic functions at
infinity with their germs in $\mathcal{M}(t)$ and using Lemma \ref{LA.1} to
polynomial (\ref{B.18}) we get that there exists $\varphi_{M}\in
\mathcal{M}(t)^{\ast}$ such that
\begin{align}
&  Q(\varphi_{M}(t),\lambda_{M}(t),u(t))=0\;\;\;\text{in }\mathcal{M}%
(t)^{\ast},\label{B.20}\\
&  \deg\varphi_{M}(t)=\max_{i=1}^{N}\frac{\deg Q_{i}(\lambda_{M}%
(t),u(t))-\deg(Q_{0}(\lambda_{M}(t),u(t))}{i}.\label{B.21}%
\end{align}

By (\ref{B.20}), similarly as in the proof of Theorem \ref{T4.7}, there exists
$\psi_{M}(t)\in\mathcal{M}(t)^{\ast}$ such that
\begin{align}
f(\varphi_{M}(t),\psi_{M}(t)) &  =\lambda_{M}(t),\label{B.22}\\
f_{y}^{\prime}(\varphi_{M}(t),\psi_{M}(t)) &  =u(t).\label{B.23}%
\end{align}
Then there exist an integer $D_{M}\in\mathbb{N}$ and functions $x_{M}%
(t),y_{M}(t)$ meromorphic at infinity which germs at infinity are equal to
$\varphi_{M}(t^{D_{M}}),\psi_{M}(t^{D_{M}}),$ respectively. Put $\Phi
_{M}(t):=(x_{M}(t),y_{M}(t)).$ By (\ref{B.22}) and the definition of
$\lambda_{M}(t)$ we have
\begin{equation}
\deg(f\circ\Phi_{M}-\lambda_{0})<0.\label{B.24'}%
\end{equation}
Hence and by (\ref{2.1}) we get
\begin{equation}
\deg\Phi_{M}(t)=D_{M}\deg\varphi_{M}(t).\label{B.25}%
\end{equation}
By (\ref{B.24'}) and (\ref{B.25}) we get $\deg\nabla f\circ\Phi_{M}=\deg
f_{y}^{\prime}\circ\Phi_{M}.$ Hence according to (\ref{B.23}) and (\ref{B.25})
we get
\[
\frac{\deg\nabla f\circ\Phi_{M}}{\deg\Phi_{M}}=\frac{\deg u(t)}{\deg
\varphi_{M}(t)}.
\]
This, together with (\ref{B.19}) and (\ref{B.21}), gives for almost all
$M\in\mathbb{N}$%
\[
\deg\Phi_{M}>0
\]
and
\[
\lim_{M\rightarrow\infty}\frac{\deg\nabla f\circ\Phi_{M}}{\deg\Phi_{M}}=0.
\]

Hence and (\ref{B.24'}), taking into account definition (\ref{2.1}), we get
$\mathcal{L}_{\infty,\lambda_{0}}(f)\leqslant0.$ Hence and by (\ref{4.30}) we
get the assertion of the theorem in case (b$_{2}$).

This ends the proof of the theorem.
\end{proof}

As\ a corollary of Theorems \ref{T4.1}, \ref{T4.7} and \ref{T4.8} we obtain a
theorem on the function $\mathbb{\mathbb{C\ni\lambda\longmapsto}}%
\mathcal{L}_{\infty,\lambda}(f)\in\mathbb{R\cup\{-\infty\}}.$

\begin{theorem}
\label{T4.9}The function $\mathbb{\mathbb{C\ni\lambda\longmapsto}}%
\mathcal{L}_{\infty,\lambda}(f)\mathbb{\ }$takes the values in $[-\infty,-1)$
for $\lambda\in\Lambda(f).$ Outside $\Lambda(f)$ this function is constant and
non-negative. Moreover,

\begin{description}
\item[(a)] if $\deg_{u}Q_{0}(\lambda,u)=0,$ then $\mathcal{L}_{\infty,\lambda
}(f)=const.>0$ \ \ for \ $\lambda\notin\Lambda(f),$

\item[(b)] if $\deg_{u}Q_{0}(\lambda,u)>0,$ then $\mathcal{L}_{\infty,\lambda
}(f)=const.=0$ \ \ for \ $\lambda\notin\Lambda(f).$
\end{description}
\end{theorem}

\begin{flushright}$\square$\end{flushright}

At the end we illustrate Theorem \ref{T4.9} with three simple examples.

\begin{example}

\begin{description}
\item[(a)] For $f(x,y):=y^{2}+x$ we have $\mathcal{L}_{\infty,\lambda
}(f)\equiv1.\label{E4.10}$

\item[(b)] For $f(x,y):=y^{n+1}+xy^{n}+y,$ $n>1,$ we have
\[
\mathcal{L}_{\infty,\lambda}(f)=\left\{
\begin{array}
[c]{ccc}%
\frac{1}{n} & \text{for} & \lambda\neq0,\\
-1-\frac{1}{n-1} & \text{for} & \lambda=0.
\end{array}
\right.
\]

\item[(c)] For $f(x,y):=y^{2}$ we have
\[
\mathcal{L}_{\infty,\lambda}(f)=\left\{
\begin{array}
[c]{ccc}%
0 & \text{for} & \lambda\neq0,\\
-\infty & \text{for} & \lambda=0.
\end{array}
\right.
\]
\end{description}
\end{example}

\begin{flushright}$\square$\end{flushright}

At the end of this section we compare the \L ojasiewicz exponent of $\nabla f
$ at infinity near a fibre and on the fibre. Precisely, we shall compare\ the
function $\mathbb{\mathbb{C\ni\lambda\longmapsto}}\mathcal{L}_{\infty,\lambda
}(f)\mathbb{\ }$to the function $\mathbb{\mathbb{C\ni\lambda\longmapsto}%
}\mathcal{L}_{\infty}(\nabla f|f^{-1}(\lambda))$ for a fixed polynomial $f$ of
the form (\ref{2.1}). For simplicity we put
\[
S_{\lambda}:=f^{-1}(\lambda)\;\;\;\;\text{for \ }\lambda\in\mathbb{\mathbb{C}%
}.
\]

\begin{theorem}
\label{T6.1} The functions $\mathbb{\mathbb{C\ni\lambda\longmapsto}%
}\mathcal{L}_{\infty,\lambda}(f)\mathbb{\ }$and $\mathbb{\mathbb{C\ni
\lambda\longmapsto}}\mathcal{L}_{\infty}(\nabla f|S_{\lambda})$ are identical
besides a finite set of points $\lambda_{0}$, in which $\mathcal{L}%
_{\infty,\lambda_{0}}(f)\mathbb{\ <}\mathcal{L}_{\infty}(\nabla f|S_{\lambda
_{0}}).$ This inequality holds if and only if when $\mathcal{L}_{\infty
,\lambda_{0}}(f)\in(-\infty,-1)$ or when $\mathcal{L}_{\infty,\lambda_{0}%
}(f)>0$ and
\begin{equation}
\max_{i=1}^{N}\frac{\deg_{u}Q_{i}(\lambda,u)}{i}>\max_{i=1}^{N}\frac{\deg
_{u}Q_{i}(\lambda_{0},u)}{i}.\label{*}%
\end{equation}
\end{theorem}

\begin{proof}
Notice first that by the form (\ref{2.1}) of $f$ and Lemma \ref{L3.1} we
easily get
\begin{equation}
\mathcal{L}_{\infty}(\nabla f|S_{\lambda})=\mathcal{L}_{\infty}(f_{y}^{\prime
}|S_{\lambda})\;\;\;\text{for }\lambda\in\mathbb{C}.\label{6.1}%
\end{equation}

Let $Q(x,\lambda,u)=Q_{0}(\lambda,u)x^{N}+\cdots+Q_{N}(\lambda,u),\;\;\;Q_{0}%
\neq0,$ be, as before, the resultant of $f(x,y)-\lambda$ and $f_{y}^{\prime
}(x,y)-u$ with respect to $y.$ Fix $\lambda_{0}\in\mathbb{C}$ and consider
four cases:

(a) $Q_{0}(\lambda_{0},0)\neq0$ and $\deg_{u}Q_{0}(\lambda,u)=0.$ Then by
Lemma \ref{LP} and (\ref{6.1}) we have%

\[
\mathcal{L}_{\infty}(\nabla f|S_{\lambda_{0}})=\left[  \max_{i=1}%
^{N}\frac{\deg_{u}Q_{i}(\lambda_{0},u)}{i}\right]  ^{-1}.
\]
Hence and by Corollary \ref{L4.3} and Theorem \ref{T4.7} we get $\mathcal{L}%
_{\infty}(\nabla f|S_{\lambda_{0}})>\mathcal{L}_{\infty,\lambda_{0}}(f)$ if
and only if (\ref{*}) holds. Obviously, the set of such points\ is finite.

(b) $Q_{0}(\lambda_{0},0)\neq0$ and $\deg_{u}Q_{0}(\lambda,u)>0.$ Then by
Lemma \ref{LP}, (\ref{6.1}) and Theorem \ref{T4.8} we have $\mathcal{L}%
_{\infty}(\nabla f|S_{\lambda_{0}})=\mathcal{L}_{\infty,\lambda_{0}}(f)=0.$

(c) There exists $r\in\{0,...,N-1\}$ such that $Q_{0}(\lambda_{0}%
,0)=\cdots=Q_{r}(\lambda_{0},0)=0$ $Q_{r+1}(\lambda_{0},0)\neq0.$ Then by
Theorem \ref{T4.1} $-\infty<\mathcal{L}_{\infty,\lambda_{0}}(f)<0. $ Hence and
by Corollary \ref{C3.4}, $\mathcal{L}_{\infty}(f-\lambda_{0},f_{y}^{\prime
})<0.$ Then by Theorems \ref{T3.3}, \ref{T3.2} and (\ref{6.1}) we obtain
\[
\mathcal{L}_{\infty,\lambda_{0}}(f)<\mathcal{L}_{\infty}(f-\lambda_{0}%
,f_{y}^{\prime})=\mathcal{L}_{\infty}(f-\lambda_{0}|S_{y})<\mathcal{L}%
_{\infty}(f_{y}^{\prime}|S_{\lambda_{0}})=\mathcal{L}_{\infty}(\nabla
f|S_{\lambda_{0}}).
\]
Obviously, the set of $\lambda_{0}\in\mathbb{C}$ for which (c) holds is
contained in $\Lambda(f)$ and thus finite.

(d) $Q_{0}(\lambda_{0},0)=\cdots=Q_{N}(\lambda_{0},0)=0.$ Then by Theorem
\ref{T4.1}, Lemma \ref{LP} and (\ref{6.1}) we have $-\infty=\mathcal{L}%
_{\infty,\lambda_{0}}(f)=\mathcal{L}_{\infty}(f_{y}^{\prime}|S_{\lambda_{0}%
})=\mathcal{L}_{\infty}(\nabla f|S_{\lambda_{0}}).$

This ends the proof.
\end{proof}

We illustrate the above theorem with the same polynomials as in Example
\ref{E4.10}.

\begin{example}
\label{E6.1}

\begin{description}
\item[(a)] For $f(x,y):=y^{2}+x$ we have $\mathcal{L}_{\infty,\lambda
}(f)=\mathcal{L}_{\infty}(\nabla f|S_{\lambda})=1$ for\ $\lambda
\in\mathbb{\mathbb{C}}.$

\item[(b)] For $f(x,y):=y^{n+1}+xy^{n}+y$ we have $\mathcal{L}_{\infty
,\lambda}(f)=\mathcal{L}_{\infty}(\nabla f|S_{\lambda})=\frac{1}{n}$ for
$\lambda\neq0$ and $-1-\frac{1}{n-1}=\mathcal{L}_{\infty,0}(f)<\mathcal{L}%
_{\infty}(\nabla f|S_{0})=0.$

\item[(c)] For $f(x,y):=y^{2}$ we have $\mathcal{L}_{\infty,\lambda
}(f)=\mathcal{L}_{\infty}(\nabla f|S_{\lambda})$ for any $\lambda
\in\mathbb{\mathbb{C}}.$
\end{description}
\end{example}

\begin{flushright}$\square$\end{flushright}

\section{Equivalence of the definitions of $\widetilde{\mathcal{L}}%
_{\infty,\lambda}(f)$ and $\mathcal{L}_{\infty,\lambda}(f)$}

In the Introduction we have defined $\widetilde{\mathcal{L}}_{\infty,\lambda
}(f)$ and $\mathcal{L}_{\infty,\lambda}(f)$ by formulas (\ref{1.1'}) and
(\ref{1.1}), respectively. We notice that the limit in (\ref{1.1'}) always
exists (it may happen to be $-\infty$) because by definition of $\mathcal{L}%
_{\infty}(\nabla f|f^{-1}(D_{\delta}))$ the function $\delta\mapsto
\mathcal{L}_{\infty}(\nabla f|f^{-1}(D_{\delta}))$ is non-increasing.

We now prove

\begin{theorem}
\label{T5.1}Let $f:\mathbb{C}^{2}\rightarrow\mathbb{\mathbb{C}}$ be a
non-constant polynomial and $\lambda_{0}\in\mathbb{C}.$ Then
\[
\widetilde{\mathcal{L}}_{\infty,\lambda_{0}}(f)=\mathcal{L}_{\infty
,\lambda_{0}}(f)
\]
holds.
\end{theorem}

\begin{proof}
Obviously
\[
\widetilde{\mathcal{L}}_{\infty,\lambda_{0}}(f)\leqslant\mathcal{L}%
_{\infty,\lambda_{0}}(f).
\]

We shall now prove the opposite inequality. Since the set $\Lambda(f)$ is
finite then there exists $D=\{\lambda\in\mathbb{\mathbb{C}}:\left|
\lambda-\lambda_{0}\right|  <\delta_{0}\}$ such that $(D\setminus\{\lambda
_{0}\})\cap\Lambda(f)=\emptyset.$ According to Theorem \ref{T4.9} we have%

\begin{equation}
\mathcal{L}_{\infty,\lambda}(f)\geqslant\mathcal{L}_{\infty,\lambda_{0}%
}(f)\text{ \ \ \ \ for \ }\lambda\in D.\label{5.1}%
\end{equation}
Take an arbitrary $0<\delta<\delta_{0}$ and put $D_{\delta}:=\{\lambda
\in\mathbb{\mathbb{C}}:\left|  \lambda-\lambda_{0}\right|  <\delta\}.$ Since
the set $f^{-1}(\overline{D}_{\delta})$ is semi-algebraic in $\mathbb{C}^{2}$
by the Curve Selection Lemma the exponent $\mathcal{L}_{\infty}(\nabla
f|f^{-1}(\overline{D}_{\delta}))$ is attained on a meromorphic curve
$\Phi_{\delta},$ $\deg\Phi_{\delta}>0,$ lying in this set (see \cite{CK3},
Proposition 1). It is easy to see that there exists $\tilde{\lambda}%
\in\overline{D}_{\delta}\subset D$ such that $\deg(f-\tilde{\lambda})\circ
\Phi_{\delta}<0.$ By definition of $\mathcal{L}_{\infty,\lambda_{0}}(f)$ and
(\ref{5.1}) we get $\mathcal{L}_{\infty}(\nabla f|f^{-1}(\overline{D}_{\delta
}))\geqslant\mathcal{L}_{\infty,\tilde{\lambda}}(f)\geqslant\mathcal{L}%
_{\infty,\lambda_{0}}(f).$ Hence
\[
\lim_{\delta\rightarrow0^{+}}\mathcal{L}_{\infty}(\nabla f|f^{-1}(D_{\delta
}))=\lim_{\delta\rightarrow0^{+}}\mathcal{L}_{\infty}(\nabla f|f^{-1}%
(\overline{D}_{\delta}))\geqslant\mathcal{L}_{\infty,\lambda_{0}}(f).
\]

This ends the proof.
\end{proof}

\section{$n$-dimensional case}

Let us start with definitions. A non-constant polynomial $f:\mathbb{C}%
^{n}\rightarrow\mathbb{\mathbb{C}},$ $n\geqslant2,$ is said to satisfy the
Malgrange condition for a value $\lambda_{0}\in\mathbb{\mathbb{C}}$ if
\begin{equation}
\exists_{\eta_{0},\delta_{0},R_{0}>0}\forall_{p\in\mathbb{C}^{n}}(\left|
p\right|  >R_{0}\;\wedge\;\left|  f(p)-\lambda_{0}\right|  <\delta
_{0}\;\Rightarrow\;\left|  p\right|  \left|  \nabla f(p)\right|  >\eta
_{0}).\label{7.1}%
\end{equation}
By $K_{\infty}(f)$ we denote the set of $\lambda\in\mathbb{\mathbb{C}}$ for
which the Malgrange condition does not hold. It is easy to check that
$\lambda\in K_{\infty}(f)$ if and only if there exists a sequence
$\{p_{k}\}\subset\mathbb{C}^{n}$ such that
\begin{equation}
\lim_{k\rightarrow\infty}\left|  p_{k}\right|  =\infty,\;\lim_{k\rightarrow
\infty}f(p_{k})=\lambda,\;\text{and }\lim_{k\rightarrow\infty}\left|
p_{k}\right|  \left|  \nabla f(p_{k})\right|  =0.\label{7.2}%
\end{equation}

It is known (see \cite{JK}, cf. \cite{S}) that
\begin{equation}
\#K_{\infty}(f)<+\infty.\label{7.3}%
\end{equation}

We now give a characterization of $K_{\infty}(f)$ in terms of the exponents
$\widetilde{\mathcal{L}}_{\infty,\lambda}(f)$ and $\mathcal{L}_{\infty
,\lambda}(f).$

\begin{theorem}
\label{T7.1}For $\lambda_{0}\in\mathbb{\mathbb{C}}$ the following conditions
are equivalent:

\begin{description}
\item[(i)] $\lambda_{0}\in K_{\infty}(f),$

\item[(ii)] $\widetilde{\mathcal{L}}_{\infty,\lambda_{0}}(f)<-1,$

\item[(iii)] $\mathcal{L}_{\infty,\lambda_{0}}(f)<-1.$
\end{description}
\end{theorem}

\begin{proof}
(iii)$\Rightarrow$(ii)$\Rightarrow$(i). Take $\lambda_{0}\notin K_{\infty
}(f).$ Then $\lambda_{0}$ satisfies (\ref{7.1}). Without loss of generality,
by (\ref{7.3}), we may assume that $\{\lambda\in\mathbb{\mathbb{C}}:\left|
\lambda-\lambda_{0}\right|  <\delta_{0}\}\cap K_{\infty}(f)=\emptyset.$ Take
$\delta,$ $0<\delta<\delta_{0},$ and put $D_{\delta}:=\{\lambda\in
\mathbb{\mathbb{C}}:\left|  \lambda-\lambda_{0}\right|  <\delta\}.$ Then by
(\ref{7.1}) we have $\mathcal{L}_{\infty}(\nabla f|f^{-1}(D_{\delta
}))\geqslant-1.$ Hence by definition (\ref{1.1'}) we get $\widetilde
{\mathcal{L}}_{\infty,\lambda_{0}}(f)\geqslant-1.$ From an obvious inequality
\begin{equation}
\mathcal{L}_{\infty,\lambda}(f)\geqslant\widetilde{\mathcal{L}}_{\infty
,\lambda}(f)\;\;\;\text{for \ }\lambda\in\mathbb{\mathbb{C}}\label{7.3'}%
\end{equation}
we also get $\mathcal{L}_{\infty,\lambda_{0}}(f)\geqslant-1.$ This gives the
required sequence of implications.

We now show the implication (i)$\Rightarrow$(iii). Let $\lambda_{0}\in
K_{\infty}(f)$ and $\{p_{k}\}\subset\mathbb{C}^{n}$ be a sequence satisfying
(\ref{7.2}). By (\ref{7.3}) there exists a closed disc $\overline
{D}:=\{\lambda\in\mathbb{\mathbb{C}}:\left|  \lambda-\lambda_{0}\right|
\leqslant\delta_{0}\}$ such that $\overline{D}\cap K_{\infty}(f)=\emptyset.$
Since $f^{-1}(\overline{D})$ is a semialgebraic set in $\mathbb{C}^{n},$ then
by the Curve Selection Lemma the exponent $\mathcal{L}_{\infty}(\nabla
f|f^{-1}(\overline{D}))$ is attained on a meromorphic curve $\Phi,$ $\deg
\Phi>0,$ lying in this set (cf. \cite{CK3}, Proposition 1). Thus there exists
a $\tilde{\lambda}\in\overline{D}$ such that $\deg(f-\tilde{\lambda})\circ
\Phi<0.$ On the other hand almost all elements of the sequence $\{p_{k}\}$ lie
in $f^{-1}(\overline{D}).$ Then (\ref{7.2}) implies $\mathcal{L}_{\infty
}(\nabla f|f^{-1}(\overline{D}))<-1.$ In consequence $\deg\nabla f\circ
\Phi/\deg\Phi=\mathcal{L}_{\infty}(\nabla f|f^{-1}(\overline{D}))<-1.$ Hence
we get $\tilde{\lambda}\in K_{\infty}(f)$ and thus $\tilde{\lambda}%
=\lambda_{0}.$ Summing up, there exists a meromorphic curve $\Phi,$ $\deg
\Phi>0,$ such that $\deg(f-\lambda_{0})\circ\Phi<0$ and $\deg\nabla f\circ
\Phi/\deg\Phi<-1.$ Then by definition (\ref{1.1}) we have
\[
\mathcal{L}_{\infty,\lambda_{0}}(f)<-1.
\]
This gives the desired implication and ends the proof.
\end{proof}

To illustrate the usefulness of $\mathcal{L}_{\infty,\lambda}(f)$ we show how
with the help of this exponent one can find the set $K_{\infty}(f)$ in an
example. We consider the Rabier's example (see \cite{RA}, Remark 9.1).

\begin{example}
Let $f^{R}:\mathbb{C}^{3}\rightarrow\mathbb{C},$ $\ f^{R}(x,y,z):=(xy-1)yz$.
Then $K_{\infty}(f^{R})=\{0\}$ and $\mathcal{L}_{\infty,0}(f^{R})=-\infty.$

Indeed, we first show $0\in K_{\infty}(f^{R})$. Taking $\tilde{\Phi
}(t):=(t,1/t,0),$ we have $\deg\tilde{\Phi}>0,$ $f^{R}\circ\tilde{\Phi
}(t)\equiv0$ and $\deg\nabla f^{R}\circ\tilde{\Phi}=-\infty.$ Hence according
to (\ref{1.1}) we get $\mathcal{L}_{\infty,0}(f^{R})=-\infty$ and thus $0\in
K_{\infty}(f^{R})$. To prove the opposite inclusion assume that there exists
$\lambda\neq0$ such that $\lambda\in K_{\infty}(f^{R})$. Then by Theorem
\ref{T7.1} $\mathcal{L}_{\infty,\lambda}(f^{R})<-1.$ Then there exists a
meromorphic curve $\Phi=(\varphi_{1},\varphi_{2},\varphi_{3})$ such that
$\deg\Phi>0$ and
\begin{align}
\deg(f^{R}-\lambda)\circ\Phi &  <0,\label{7.9}\\
\deg\nabla f^{R}\circ\Phi &  <-\deg\Phi.\label{7.10}%
\end{align}
From (\ref{7.9})
\begin{equation}
\deg((\varphi_{1}\varphi_{2}-1)\varphi_{2}\varphi_{3})=0,\label{7.11}%
\end{equation}
whereas from (\ref{7.10}) we get $\deg f_{z}^{\prime}\circ\Phi<-\deg\Phi$ and
thus
\begin{equation}
\deg(\varphi_{1}\varphi_{2}-1)\varphi_{2})<-\deg\Phi.\label{7.12}%
\end{equation}
By (\ref{7.11}) and (\ref{7.12}) we get $-\deg\varphi_{3}<-\deg\Phi,$ which is impossible.
\end{example}

\begin{flushright}$\square$\end{flushright}

Using $\mathcal{L}_{\infty,\lambda}(f)$ we shall prove one more theorem

\begin{theorem}
\label{T6.1.2}Let $f:\mathbb{C}^{n}\rightarrow\mathbb{\mathbb{C}}$,
$n\geqslant2,$ be a non-constant polynomial. If $\mathcal{L}_{\infty}(\nabla
f)\leqslant-1,$ then there exists $\lambda_{0}\in\mathbb{C}$ such that
\begin{equation}
\mathcal{L}_{\infty}(\nabla f)=\mathcal{L}_{\infty,\lambda_{0}}(f).\label{A.5}%
\end{equation}
\end{theorem}

\begin{proof}
Let $\Phi,$ $\deg\Phi>0,$ be a meromorphic curve on which the \L ojasiewicz
exponent $\mathcal{L}_{\infty}(\nabla f)$ is attained. Then
\begin{equation}
\mathcal{L}_{\infty}(\nabla f)=\frac{\deg\nabla f\circ\Phi}{\deg\Phi
}.\label{A.6}%
\end{equation}
We shall show
\begin{equation}
\deg f\circ\Phi\leqslant0.\label{A.7}%
\end{equation}
Indeed, it suffices to consider the case $\deg f\circ\Phi\neq0.$ Then
\[
\frac{\deg f\circ\Phi}{\deg\Phi}\leqslant\frac{\deg\nabla f\circ\Phi}{\deg
\Phi}+1=\mathcal{L}_{\infty}(\nabla f)+1\leqslant0,
\]
which gives (\ref{A.7}).

Inequality (\ref{A.7}) implies that there exists $\lambda_{0}\in\mathbb{C}$
such that
\begin{equation}
\deg(f-\lambda_{0})\circ\Phi<0.\label{A.8}%
\end{equation}
Then by (\ref{A.6}), (\ref{A.8}) and (\ref{1.1}) we get
\[
\mathcal{L}_{\infty,\lambda_{0}}(f)\leqslant\mathcal{L}_{\infty}(\nabla f).
\]

The opposite inequality is obvious.

This ends the proof.
\end{proof}

Directly from the above theorem we obtain

\begin{theorem}
Let $f:\mathbb{C}^{n}\rightarrow\mathbb{\mathbb{C}}$, $n\geqslant2,$ be a
non-constant polynomial. The following conditions are equivalent:

\begin{description}
\item[(i)] $K_{\infty}(f)\neq\emptyset,$

\item[(ii)] $\mathcal{L}_{\infty}(\nabla f)<-1.$
\end{description}
\end{theorem}

\begin{proof}
(i)$\Rightarrow$(ii). Take $\lambda_{0}\in K_{\infty}(f).$ Then by Theorem
\ref{T7.1} we have $\mathcal{L}_{\infty,\lambda_{0}}(f)<-1.$ Then
$\mathcal{L}_{\infty}(\nabla f)<-1.$

(ii)$\Rightarrow$(i). By Theorem \ref{T6.1.2} there exists $\lambda_{0}%
\in\mathbb{C}$ such that $\mathcal{L}_{\infty,\lambda_{0}}(f)=\mathcal{L}%
_{\infty}(\nabla f)<-1.$ Hence by Theorem \ref{T7.1} $\lambda_{0}\in
K_{\infty}(f).$

This ends the proof.
\end{proof}

\begin{theorem}
\label{T7.2}We have
\begin{equation}
\widetilde{\mathcal{L}}_{\infty,\lambda}(f)=\mathcal{L}_{\infty,\lambda
}(f)\;\;\;\text{for \ }\lambda\in K_{\infty}(f).\label{7.5}%
\end{equation}
\end{theorem}

\begin{proof}
Take any $\lambda_{0}\in K_{\infty}(f).$ Since the inequality $\mathcal{L}%
_{\infty,\lambda_{0}}(f)\geqslant\widetilde{\mathcal{L}}_{\infty,\lambda_{0}%
}(f)$ is obvious it suffices to show
\begin{equation}
\mathcal{L}_{\infty,\lambda_{0}}(f)\leqslant\widetilde{\mathcal{L}}%
_{\infty,\lambda_{0}}(f).\label{7.6}%
\end{equation}

Similarly as in the second part of the proof of Theorem \ref{T7.1} we take a
closed disc $\overline{D}_{\delta}:=\{\lambda\in\mathbb{\mathbb{C}}:\left|
\lambda-\lambda_{0}\right|  \leqslant\delta\}$ such that $\overline{D}%
_{\delta}\cap K_{\infty}(f)=\emptyset.$ Then analogously as previously we show
that there exists a meromorphic curve $\Phi_{\delta}$ lying in $f^{-1}%
(\overline{D}_{\delta})$ such that,
\begin{equation}
\deg\Phi_{\delta}>0,\;\deg(f-\lambda_{0})\circ\Phi_{\delta}<0,\;\mathcal{L}%
_{\infty}(\nabla f|f^{-1}(\overline{D}_{\delta}))=\frac{\deg\nabla f\circ
\Phi_{\delta}}{\deg\Phi_{\delta}}.\label{7.7}%
\end{equation}
Hence and from (\ref{1.1}) we get
\begin{equation}
\mathcal{L}_{\infty,\lambda_{0}}(f)\leqslant\mathcal{L}_{\infty}(\nabla
f|f^{-1}(\overline{D}_{\delta}))\label{7.8}%
\end{equation}
Then taking $\delta\rightarrow0^{+}$ we get (\ref{7.6}).

This ends the proof.
\end{proof}

We give now a complement of Theorem \ref{T7.2}, which follows immediately from
Theorem \ref{T7.1}.

\begin{proposition}
\label{P6.2.1}If $\lambda_{0}\notin K_{\infty}(f)$ and there exists a
meromorphic curve such that $\deg(f-\lambda_{0})\circ\Phi<0$ and $\deg\nabla
f\circ\Phi=-\deg\Phi$ then
\[
\mathcal{L}_{\infty,\lambda_{0}}(f)=\widetilde{\mathcal{L}}_{\infty
,\lambda_{0}}(f)=-1.
\]
\end{proposition}

\begin{flushright}$\square$\end{flushright}

Before we pose questions concerning $\widetilde{\mathcal{L}}_{\infty,\lambda
}(f)$ and $\mathcal{L}_{\infty,\lambda}(f)$ for $n>2$ we return to the
Rabier's example.

\begin{example}
\label{ER}For the Rabier's polynomial we have
\begin{equation}
\mathcal{L}_{\infty,\lambda}(f^{R})=\widetilde{\mathcal{L}}_{\infty,\lambda
}(f^{R})=\left\{
\begin{array}
[c]{ccc}%
-\infty & \text{for} & \lambda=0,\\
-1 & \text{for} & \lambda\neq0.
\end{array}
\right. \label{y}%
\end{equation}
Indeed, for any $\lambda\in\mathbb{C}$ and $\Phi_{\lambda}%
(t):=(t,1/2t,-4\lambda t)$ we have $f^{R}\circ\Phi_{\lambda}(t)\equiv\lambda$
and $\deg\nabla f\circ\Phi_{\lambda}=-\deg\Phi_{\lambda}.$ Hence and from
Proposition \ref{P6.2.1} we get (\ref{y}).
\end{example}

\begin{flushright}$\square$\end{flushright}

This example gives rise to the question: can one generalize Theorems
\ref{T5.1} and \ref{T4.9} to the $n$ dimensional case? Precisely

\begin{problem}
\label{P1}Does the equality
\[
\mathcal{L}_{\infty,\lambda}(f)=\widetilde{\mathcal{L}}_{\infty,\lambda}(f)
\]
hold for any non-constant polynomial $f:\mathbb{C}^{n}\rightarrow
\mathbb{\mathbb{C}},$ $n>2,$ and $\lambda\in\mathbb{C}?$
\end{problem}

\begin{flushright}$\square$\end{flushright}

\begin{problem}
\label{P2}Does for a non-constant polynomial $f:\mathbb{C}^{n}\rightarrow
\mathbb{\mathbb{C}},$ $n>2,$ there exists a number $c_{f}\in\lbrack-1,\infty)$
such that
\[
\mathcal{L}_{\infty,\lambda}(f)=\widetilde{\mathcal{L}}_{\infty,\lambda
}(f)=c_{f}\text{ \ \ \ for any }\lambda\notin K_{\infty}(f)?
\]
\end{problem}

\begin{flushright}$\square$\end{flushright}

\begin{remark}
It is known for $n=2$ (see i.e. \cite{S}) that $\Lambda(f)=K_{\infty}(f)$.
Then the answers to both Problems \ref{P1} and \ref{P2} are positive.
Moreover, by Theorem \ref{T4.9} the constant $c_{f}$ belongs to the interval
$[0,\infty).$

For $n>2$ we have only the inclusion $\Lambda(f)\subset K_{\infty}(f)$ (see
\cite{PA}). In general the equality does not hold. One can show that for the
polynomial $f^{PZ}(x,y,z):=x-3x^{5}y^{2}+2x^{7}y^{3}+yz$ (see \cite{PZ}) we
have
\[
\emptyset=\Lambda(f^{PZ})\varsubsetneq K_{\infty}(f^{PZ})\neq\emptyset.
\]
Therefore by (\ref{7.3}) it seems to be more natural in the case $n>2$ to
consider the dependence of the exponents $\widetilde{\mathcal{L}}%
_{\infty,\lambda}(f)$ and $\mathcal{L}_{\infty,\lambda}(f)$ on $K_{\infty}(f)$
and not on $\Lambda(f).$
\end{remark}

\begin{flushright}$\square$\end{flushright}

Let us introduce one more definition. A non-constant polynomial $f:\mathbb{C}%
^{n}\rightarrow\mathbb{\mathbb{C}},$ $n\geqslant2,$ is said to satisfy the
Fedorjuk condition for a value $\lambda_{0}\in\mathbb{\mathbb{C}}$ if
\begin{equation}
\exists_{\eta_{0},\delta_{0},R_{0}>0}\forall_{p\in\mathbb{C}^{n}}(\left|
p\right|  >R_{0}\;\wedge\;\left|  f(p)-\lambda_{0}\right|  <\delta
_{0}\;\Rightarrow\;|\nabla f(p)|>\eta_{0}).\label{6.18}%
\end{equation}
By $\widetilde{K}_{\infty}(f)$ we denote the set of $\lambda\in
\mathbb{\mathbb{C}}$ for which the Fedorjuk condition does not hold. It is
easy to check that $\lambda\in\widetilde{K}_{\infty}(f)$ if and only if there
exists a sequence $\{p_{k}\}\subset\mathbb{C}^{n}$ such that
\begin{equation}
\lim_{k\rightarrow\infty}\left|  p_{k}\right|  =\infty,\;\lim_{k\rightarrow
\infty}f(p_{k})=\lambda,\;\text{and }\lim_{k\rightarrow\infty}\left|  \nabla
f(p_{k})\right|  =0.\label{6.19}%
\end{equation}

It is known (see \cite{S}) that $\widetilde{K}_{\infty}(f)$ is algebraic. Then
we have two excluded possibilities
\begin{equation}
\#\widetilde{K}_{\infty}(f)<+\infty\;\;\text{or}\;\;\widetilde{K}_{\infty
}(f)=\mathbb{\mathbb{C}}.
\end{equation}

Analogously to Theorem \ref{T7.1} we show

\begin{theorem}
\label{T6.11}If $\widetilde{K}_{\infty}(f)\not =\mathbb{\mathbb{C}}$ then for
$\lambda_{0}\in\mathbb{\mathbb{C}}$ the following conditions are equivalent:

\begin{description}
\item[(i)] $\lambda_{0}\in\widetilde{K}_{\infty}(f),$

\item[(ii)] $\widetilde{\mathcal{L}}_{\infty,\lambda_{0}}(f)<0,$

\item[(iii)] $\mathcal{L}_{\infty,\lambda_{0}}(f)<0.$
\end{description}
\end{theorem}

\begin{flushright}$\square$\end{flushright}

\begin{remark}
The inclusion $K_{\infty}(f)\subset\widetilde{K}_{\infty}(f)$ is obvious. For
$n=2$ we have $K_{\infty}(f)=\widetilde{K}_{\infty}(f)$ and thus
$\widetilde{K}_{\infty}(f)$ is finite. For $n>2$ the equality $K_{\infty
}(f)=\widetilde{K}_{\infty}(f)$ does not have to hold. Namely, from Example
\ref{ER} and Theorem \ref{T6.11} we have $\widetilde{K}_{\infty}%
(f^{R})=\mathbb{C}$ whereas $K_{\infty}(f^{R})$ is finite.
\end{remark}

\begin{flushright}$\square$\end{flushright}

\begin{problem}
Does the equality $\widetilde{K}_{\infty}(f)=\mathbb{\mathbb{C}}$ imply that
$\mathcal{L}_{\infty,\lambda}(f)<0$ for each $\lambda\in\mathbb{\mathbb{C}}$ ?
\end{problem}

\begin{flushright}$\square$\end{flushright}

\end{document}